\documentclass[thmsa,onecolumn,10pt]{article}
\usepackage{amssymb}
\usepackage{graphicx}
\usepackage{amsmath}

\newtheorem{theorem}{Theorem}

\newtheorem{corollary}[theorem]{Corollary}

\newtheorem{lemma}[theorem]{Lemma}

\newtheorem{proposition}[theorem]{Proposition}

\newenvironment{proof}[1][Proof]{\textbf{#1.} }{\ \rule{0.5em}{0.5em}}

\begin{document}

\title{\textsc{High-frequency asymptotics for subordinated isotropic fields
on an Abelian compact group\thanks{%
D.M. is grateful to Paolo Baldi for many useful discussions.}}}
\author{\textsc{Domenico MARINUCCI}\thanks{%
Department of Mathematics, University of Rome \textquotedblleft Tor
Vergata\textquotedblright . E-mail: \texttt{marinucc@mat.uniroma2.it}}\ \ \
and\ \textsc{Giovanni PECCATI}\thanks{%
Laboratoire de Statistique Th\'{e}orique et Appliqu\'{e}e de l'Universit\'{e}
Paris VI. E-mail: \texttt{giovanni.peccati@gmail.com}}}
\date{\textsc{July 3, 2006}}
\maketitle

\begin{abstract}
Let $\widetilde{T}\left( g\right) $ be a random field indexed by an Abelian
compact group $G$, and suppose that $\widetilde{T}$ has the form $\widetilde{%
T}\left( g\right) =F\left( T\left( g\right) \right) $, where $T$ is Gaussian
and isotropic. The aim of this paper is to establish high-frequency central
limit theorems for the Fourier coefficients associated to $\widetilde{T}$.
The proofs of our main results involve recently established criteria for the
weak convergence of multiple Wiener-It\^{o} integrals. Our research is
motivated by physical applications, mainly related to the probabilistic
modelization of the Cosmic Microwave Background radiation. In this
connection, the case of the $n$-dimensional torus is analyzed in detail.


\textbf{AMS\ classification -- }Primary 60B15; Secondary 60F05, 60G60

\textbf{Running Title - } \textsc{asymptotics for random fields on an
Abelian group}
\end{abstract}

\section{Introduction \label{S : Intro}}

Let $G$ be a connected compact Abelian group. The aim of this paper is to
establish central limit theorems (CLTs) for the Fourier coefficients
associated to a random field indexed by $G$, and subordinated to some
real-valued isotropic Gaussian field $T=\left\{ T\left( g\right) :g\in
G\right\} $. By \textit{isotropic} we mean that, for every $p\geq 1$ and
every $h,g_{1},...,g_{p}\in G$,%
\begin{equation}
\left\{ T\left( hg_{1}\right) ,...,T\left( hg_{p}\right) \right\} \overset{%
law}{=}\left\{ T\left( g_{1}\right) ,...,T\left( g_{p}\right) \right\} \text{%
,}  \label{isotrop}
\end{equation}%
i.e. the finite-dimensional distributions of the \textquotedblleft
translated\textquotedblright\ process $g\mapsto T\left( hg\right) $ coincide
with those of $T$, for every $h\in G$. As a consequence of the Peter-Weyl
theorem (see e.g. \cite{DUistK}), the Gaussian field $T$ always admits the
expansion%
\begin{equation}
T\left( g\right) =\sum_{\pi \in \hat{G}}a_{\pi }\chi _{\pi }\left( g\right)
\text{, \ \ }g\in G\text{,}  \label{start}
\end{equation}%
where $\hat{G}$ is the collection of the irreducible unitary representations
of $G$ (that is, $\hat{G}$ is the \textit{dual }of $G$ -- see e.g. \cite%
{Rudin}), $\chi _{\pi }$ is the \textit{character }associated to a given $%
\pi \in \hat{G}$, and%
\begin{equation}
a_{\pi }\triangleq \int_{G}T\left( g\right) \chi _{\pi }\left( g^{-1}\right)
dg  \label{Fcooeff}
\end{equation}%
with $dg$ indicating the Haar measure (a more detailed discussion of the
properties of the expansion (\ref{Fcooeff}) is deferred to the next
section). Now consider a real-valued $F\in L^{2}\left( \mathbb{R},\exp
\left( -x^{2}/2\right) dx\right) $, and define the \textit{subordinated
field }$F\left[ T\right] $ as
\begin{equation}
F\left[ T\right] \left( g\right) \triangleq F\left( T\left( g\right) \right)
\text{, \ \ }\forall g\in G\text{.}  \label{effeti}
\end{equation}%
Plainly, for a non-linear transformation $F$ the field $F\left[ T\right] $
is in general not Gaussian. However, since $T$ is isotropic $F\left[ T\right]
$ is isotropic, and the Peter-Weyl Theorem yields again the spectral
expansion
\begin{equation}
F\left[ T\right] \left( g\right) =\sum_{\pi \in \hat{G}}\widetilde{a}_{\pi
}\left( F\right) \chi _{\pi }\left( g\right) \text{, \ \ }g\in G\text{,}
\label{subFourier}
\end{equation}%
where
\begin{equation}
\widetilde{a}_{\pi }\left( F\right) \triangleq \int_{G}F\left[ T\right]
\left( g\right) \chi _{\pi }\left( g^{-1}\right) dg\text{.}  \label{FCoef}
\end{equation}%
Our aim in this paper is to investigate the asymptotic behavior of the
complex-valued variable $\widetilde{a}_{\pi }\left( F\right) $, whenever the
dual set $\hat{G}$ is infinite. More precisely, we shall establish
sufficient (and in many cases, also necessary) conditions for the following
CLT to hold:%
\begin{equation}
\mathbb{E}\left[ \left\vert \widetilde{a}_{\pi }\left( F\right) \right\vert
^{2}\right] ^{-\frac{1}{2}}\widetilde{a}_{\pi }\left( F\right) \underset{%
\left\{ \pi \right\} }{\overset{law}{\rightarrow }}N+iN^{\prime },
\label{CLT}
\end{equation}%
where $N$ and $N^{\prime }$ are two independent centered Gaussian random
variables with common variance equal to $1/2$. In (\ref{CLT}), and for the
rest paper, the subscript $\left\{ \pi \right\} $ means that $\left\{ \pi
\right\} =\left\{ \pi _{l}:l=1,2,...\right\} $ is an infinite sequence of
elements of $\hat{G}$, and that the limit is taken as $l\rightarrow +\infty $%
. A central limit result such as (\ref{CLT}) is called a \textit{%
high-frequency central limit theorem}, in analogy with the case of $G$ being
a the $n$-dimensional torus $\mathbb{R}^{n}/\left( 2\pi \mathbb{Z}\right)
^{n}$. Indeed, in this case one has that: (i) $\hat{G}$ can be identified
with the class of complex-valued mappings of the type $\mathbf{\vartheta }%
\mapsto \exp \left( i\mathbf{k}^{\prime }\mathbf{\vartheta }\right) $, where
$\mathbf{k}\in \mathbb{Z}^{n}$ and $\mathbf{\vartheta }\in (0,2\pi ]^{n}$,
(ii) the class $\left\{ \widetilde{a}_{\pi }\left( F\right) :\pi \in \hat{G}%
\right\} $ reduces to the collection of the coefficients $\left\{ \widetilde{%
a}_{\mathbf{k}}\left( F\right) :\mathbf{k}\in \mathbb{Z}^{n}\right\} $
appearing in the usual Fourier expansion $F\left[ T\right] \left( \mathbf{%
\vartheta }\right) =\sum_{\mathbf{k\in }\mathbb{Z}^{n}}\widetilde{a}_{%
\mathbf{k}}\left( F\right) \exp \left( i\mathbf{k}^{\prime }\mathbf{%
\vartheta }\right) $, and (iii) the subscript $\left\{ \pi \right\} $ in (%
\ref{CLT}) may be replaced by the condition $\left\Vert \mathbf{k}%
\right\Vert _{\mathbb{Z}^{n}}\rightarrow +\infty $, where $\left\Vert
\mathbf{\cdot }\right\Vert _{\mathbb{Z}^{n}}$ stands for the Euclidean norm.

\bigskip

Our work is strongly motivated by physical applications; indeed,
nonlinear transformations of Gaussian random fields emerge quite
naturally in a variety of physical models. A particularly active
area has recently been related to theoretical Cosmology, and more
precisely, to so-called inflationary models aimed at the
investigation of the dynamics of the gravitational potential
around the Big Bang (see for instance \cite{dodelson} and
\cite{PeebBook}). In this area, the aim is the understanding of
the primordial fluctuations which have provided the seeds for the
large scale structure of the Universe as it is currently observed,
i.e., the formation of structures such as clusters of galaxies,
filaments, walls and all those inhomogeneities which have made our
own existence possible. The currently favored scenario suggests
that the primordial seeds for these inhomogeneities have actually
been provided by quantum fluctuations in the gravitational
potential, which have then been \textquotedblleft
freezed\textquotedblright\ as large scale fluctuations when the
Universe experienced a phase of superluminal expansion known as
inflation. In these models, the primordial gravitational potential
is represented as a Gaussian field undergoing a small nonlinear
perturbation, the simplest example being provided by the so-called
Bardeen's potential
\begin{equation}
\widetilde{\Phi }(\mathbf{\vartheta })=\Phi (\mathbf{\vartheta }%
)+f_{NL}(\Phi ^{2}(\mathbf{\vartheta })-\mathbb{E}\Phi ^{2}(\mathbf{%
\vartheta }))\text{ , }\mathbf{\vartheta }\in \Theta \text{ ,}
\label{bardeen}
\end{equation}%
where $\Phi (\mathbf{\vartheta })$ denotes a zero-mean, isotropic Gaussian
random field, with parameter space $\Theta $; the nonlinearity parameter $%
f_{NL}$ can be usually described explicitly in terms of fundamental physical
constants. There is now an enormously vast physical literature on these
Gaussian subordinated fields, see for instance \cite{Bart2}, \cite{Maldac};
a recent and comprehensive survey is in \cite{Bart}. The topological
structure of $\Theta $ can vary across different physical models and it is
not unusual to assume that $\mathbf{\vartheta }$ belongs to the
three-dimensional torus $\mathbb{R}^{3}/(2\pi \mathbb{Z)}^{3}$ (see for
instance \cite{chiang}, \cite{coles}).

Very recently it has become possible to place tight observational
constraints on the predictions of inflationary models, by means of
observations on the Cosmic Microwave Background radiation (CMB). CMB\ can be
viewed as a snapshot of the Universe at the time of recombination, i.e.
\textquotedblleft soon after\textquotedblright\ the Big Bang (see again \cite%
{dodelson} for more detailed statements). It is directly related to the
primordial gravitational potential, by means of a filtering equation known
as the radiation transfer function. In the last few years huge satellite
experiments by NASA and ESA have reached the level of resolution where
models like (\ref{bardeen}) can be tested on the observations. A vast
literature has focussed on such testing procedures (for instance \cite%
{Cabella}, \cite{Kogo}, \cite{Marinucci}). An important feature of these
procedures is their asymptotic behavior; in this framework, asymptotic is
meant in the so-called high resolution sense, i.e. with respect to
observations corresponding to frequencies which become higher and higher as
the resolution of the experiment improves. On these components much effort
for physical investigation is focussing, and it is therefore of fundamental
importance to understand what is the high-frequency behavior of Gaussian
subordinated fields (see also \cite{Baldi} for other statistical
motivations). The present paper is a contribution in this direction; in
future work we shall address related issues for random fields defined on
homogenous spaces of non-Abelian groups, primarily the rotation group $SO(n)$%
, see \cite{MaPePr}.

\bigskip

The proofs of our main results rely on the classic representation
of the function $F(.)$ in (\ref{effeti}) as an infinite series of
Hermite polynomials, and on recently established criteria for the
weak convergence
of multiple Wiener-It\^{o} integrals -- as proved in \cite{NuPe} and \cite%
{PT}. Our methodology, which involves the explicit computation of the norms
associated to contraction operators, should be compared with the classic
\textquotedblleft method of diagrams\textquotedblright\ (see e.g. \cite%
{BreMaj}, \cite{GirSur} and \cite{Surg}).

\bigskip

The plan of the paper is as follows: in Section 2 we introduce our general
setting and we review some background material on random fields on groups.
Section 3 is devoted to the statements of our main results, whose proofs are
collected in Section 5, which builds upon background material on weak
convergence of multiple stochastic integrals which is collected in Section
4. Section 6 addresses some joint convergence issues, whereas Section 7 is
devoted to the analysis of general, square integrable transforms. Finally
Section 8 specializes our results to the case of the $n$-dimensional torus,
discussing the possible fulfillment of our necessary and sufficient
conditions for the CLT by physically motivated models.

\section{General setting\label{S : Gen7}}

Given $z\in \mathbb{C}$, $\Re \left( z\right) $ and $\Im \left( z\right) $
stand, respectively, for the real and imaginary part of $z$. Let $\left( G,%
\mathbb{G}\right) $ be a topological compact connected Abelian group, where $%
\mathbb{G}$ is a topology with a countable basis. As in formula (\ref{start}%
), we shall denote by $\hat{G}$ the \textit{dual }of $G$, i.e. $\hat{G}$ is
the collection of all the equivalence classes of the unitary irreducible
representations of $G$. The elements of $\hat{G}$ are noted $\pi $, $\sigma $%
, ...; the associated characters are written $\chi _{\pi }$, $\chi _{\sigma
} $, and so on. It is well known that, since $G$ is Abelian, every
irreducible representation of $G$ has dimension one. Moreover, since $G$ is
second countable (and therefore metrizable), $\hat{G}$ is at most countable.
Recall also that $\hat{G}$ is itself an Abelian group (which in general
fails to be compact), under the commutative group operation
\begin{equation}
\left( \pi ,\sigma \right) \mapsto \pi \sigma \triangleq \pi \otimes \sigma
\text{,}  \label{AbProd}
\end{equation}%
where $\otimes $ indicates the tensor product between representations. The
identity element of $\hat{G}$ is $\pi _{0}$, i.e. the trivial
representation, and $\pi ^{-1}=\overline{\pi }$, where $\overline{\pi }$
indicates complex conjugation. By using this notation, $\forall \sigma ,\pi
\in \hat{G}$ one has the obvious relations%
\begin{equation}
\chi _{\pi }\chi _{\sigma }=\chi _{\pi \sigma }\text{ \ \ and \ }\overline{%
\chi _{\pi }}=\chi _{\pi ^{-1}}=\chi _{\overline{\pi }}\text{;}
\label{AbCar}
\end{equation}%
moreover, by connectedness, $\chi _{\pi }$ is real-valued if, and only if, $%
\pi =\pi _{0}$. Observe that, since every $\pi \in \hat{G}$ has dimension
one, the distinction between $\pi $ and $\chi _{\pi }$ is immaterial (see
e.g. \cite[Corollary 4.1.2]{DUistK}). However, part of the results of this
paper can be extended to the case of a non-commutative compact group (as the
group of rotations $SO\left( 3\right) $ -- see e.g. \cite{MaPePr}) and, to
facilitate the connection between the two frameworks, we choose to adopt
this slightly redundant notation throughout Sections \ref{S : Gen7} to \ref%
{S : genCLT}.\ We note $dg$ the unique Haar measure with mass $1$ associated
to $G$, and write $L^{2}\left( G\right) =L^{2}\left( G,dg\right) $ to
indicate the space of complex-valued functions on $G$ that are
square-integrable with respect to $dg$. Since $G$ is Abelian, the class $%
\left\{ \chi _{\pi }:\pi \in \hat{G}\right\} $ is an orthonormal basis of $%
L^{2}\left( G\right) $. In what follows, $G$ will always indicate a
topological compact group such that the cardinality of $\hat{G}$ is
infinite. The reader is referred e.g. to \cite{Diaconis}, \cite[Chapter IV]%
{DUistK} or \cite{JamesLiebeck}, for every unexplained notion or result
concerning group representations.

\bigskip

We now consider a centered real-valued Gaussian random field $T=\left\{
T\left( g\right) :g\in G\right\} $ which is isotropic in the sense of
relation (\ref{isotrop}), and we shall assume for simplicity that $\mathbb{E}%
\left[ T\left( g\right) ^{2}\right] =1$. As discussed in the introduction,
the Peter-Weyl theorem implies that the spectral expansion (\ref{start})
holds, where the convergence takes place in $L^{2}\left( \Omega \times G,%
\mathbb{P\times }dg\right) .$ Note also that, for every fixed $g\in G$, the
RHS of (\ref{start}) converges in $L^{2}\left( \mathbb{P}\right) $ (see e.g.
\cite{PePy} for general results concerning decompositions of isotropic
fields).

\bigskip

Due to the isotropic and Gaussian assumptions, the class of random variables
$\left\{ a_{\pi }:\pi \in \hat{G}\right\} $ appearing in (\ref{Fcooeff}) has
a special structure (compare \cite{BaMa}). This point is summarized in the
following Lemma.

\begin{lemma}
\label{L : coeffTorus}The family $\left\{ a_{\pi }:\pi \in \hat{G}\right\} $
is composed of complex-valued Gaussian random variables such that

\begin{enumerate}
\item $a_{\pi }=\overline{a}_{\pi ^{-1}}$ for every $\pi \in \hat{G}$ (in
particular, $a_{\pi _{0}}$ is real);

\item For any $\pi ,\sigma \in \hat{G}$ such that $\pi \mathbf{\notin }%
\left\{ \sigma ,\sigma ^{-1}\right\} $, the coefficients $a_{\pi }$ and $%
a_{\sigma }$ are independent;

\item For every $\pi \neq \pi _{0}$, the random variables $\Re \left( a_{\pi
}\right) $ and $\Im \left( a_{\pi }\right) $ are Gaussian, independent,
centered and identically distributed (in particular, $\mathbb{E}\Re \left(
a_{\pi }\right) ^{2}=\mathbb{E}\Im \left( a_{\pi }\right) ^{2}$);

\item By noting
\begin{equation}
C_{\pi }\triangleq \mathbb{E}\left\vert a_{\pi }\right\vert ^{2}=2\mathbb{E}%
\left( \Re \left( a_{\pi }\right) ^{2}\right) =2\mathbb{E}\left( \Im \left(
a_{\pi }\right) ^{2}\right) ,\text{ \ \ }\pi \in \hat{G},  \label{vark}
\end{equation}%
one has $C_{\pi }=C_{\pi ^{-1}}$ and $\sum_{\pi \in \hat{G}}C_{\pi }<+\infty
.$
\end{enumerate}
\end{lemma}

\begin{proof}
Point 1 is a consequence of (\ref{Fcooeff}). The isotropic assumption
implies that $\forall \pi ,\sigma \in \hat{G}$ such that $\pi \neq \sigma $,
$\mathbb{E}\left[ a_{\pi }\overline{a_{\sigma }}\right] =0$. It follows
that, if $\pi \notin \left\{ \sigma ,\sigma ^{-1}\right\} $, $0=\mathbb{E}%
\left[ a_{\pi }\overline{a_{\sigma }}\right] =\mathbb{E}\left[ a_{\pi }%
\overline{a_{\sigma ^{-1}}}\right] =\mathbb{E}\left[ a_{\pi }a_{\sigma }%
\right] $, thus giving Point 2. Now fix $\pi \neq \pi _{0}$. Point 2,
implies that%
\begin{eqnarray*}
0 &=&\mathbb{E}\left[ a_{\pi }\overline{a_{\pi ^{-1}}}\right] =\mathbb{E}%
\left[ a_{\pi }a_{\pi }\right] \\
&=&\mathbb{E}\left( \Re \left( a_{\pi }\right) ^{2}\right) -\mathbb{E}\left(
\Im \left( a_{\pi }\right) ^{2}\right) +2i\mathbb{E}\left( \Re \left( a_{\pi
}\right) \Im \left( a_{\pi }\right) \right) \text{,}
\end{eqnarray*}%
giving immediately Point 3. Point 4 follows by combining Point 1 and Point 3.
\end{proof}

\bigskip

\textbf{Remarks -- }(a) The law of a collection of random variables $\left\{
a_{\pi }:\pi \in \hat{G}\right\} \in \mathbb{C}^{\hat{G}}$ satisfying Points
1-3 of Lemma \ref{L : coeffTorus} is completely determined by the
coefficients $C_{\pi }$ defined in (\ref{vark}).

(b) Given a collection $\left\{ a_{\pi }:\pi \in \hat{G}\right\} \in \mathbb{%
C}^{\hat{G}}$, satisfying Points 1-3 of Lemma \ref{L : coeffTorus} and such
that $\sum_{\pi \in \hat{G}}C_{\pi }<+\infty $, we may always define a
real-valued Gaussian isotropic random field $\overline{T}$ by setting $%
\overline{T}\left( g\right) $ $=$ $\sum_{\pi }a_{\pi }\chi _{\pi }\left(
g\right) $.

\bigskip

Throughout the paper, we will systematically work under the following
assumption.

\bigskip

\textbf{Assumption I }-- Let $\left\{ a_{\pi }:\pi \in \hat{G}\right\} $ be
the Fourier coefficients defined in formula (\ref{Fcooeff}), and let $%
\{C_{\pi }:\pi \in \hat{G}\}$ be given by (\ref{vark}). Then, $C_{\pi }>0$
for every $\pi \in \hat{G}$ (or, equivalently, $a_{\pi }\neq 0$, a.s.-$%
\mathbb{P}$, for every $\pi \in \hat{G}$).

\bigskip

Assumption I is a mild regularity condition on the behavior of the spectral
density of $T$. Basically, it ensures that every field of the type $g\mapsto
F\left( T\left( g\right) \right) $, where $F$ is a polynomial, admits an
expansion of the type (\ref{subFourier}) such that $\widetilde{a}_{\pi
}\left( F\right) \neq 0$ for every $\pi \in \hat{G}$, and therefore that the
asymptotic behavior of the $\widetilde{a}_{\pi }\left( F\right) $'s is not
trivial at the limit. Observe that the results of this paper extend easily
to the case of a Gaussian field $T$, such that $a_{\pi }\neq 0$ for
infinitely many $\pi $'s (at the cost of some heavier notation).

\bigskip

We now note $L_{0}^{2}\left( \mathbb{R},\exp \left( -x^{2}/2\right)
dx\right) $ the class of real-valued functions on $\mathbb{R}$, such that $%
\int_{\mathbb{R}}F\left( x\right) e^{-x^{2}/2}dx$ $=$ $0$. For a fixed $F\in
L_{0}^{2}\left( \mathbb{R},\exp \left( -x^{2}/2\right) dx\right) $, we
define the (centered) subordinated field\textit{\ }$F\left[ T\right] $ as in
(\ref{effeti}). As indicated in the introduction, $F\left[ T\right] $ is
isotropic and admits the spectral representation (\ref{subFourier}), where
the convergence of the series takes place in $L^{2}\left( \Omega \times G,%
\mathbb{P}\times dg\right) $, and, for every fixed $g\in G$, in $L^{2}\left(
\mathbb{P}\right) $. It is evident that the coefficients $\widetilde{a}_{\pi
}\left( F\right) $, $\pi \in \hat{G}$, defined in (\ref{FCoef}) are
complex-valued, centered and square integrable random variables for every $%
\pi $, and also that $\Im \left( \widetilde{a}_{\pi _{0}}\left( F\right)
\right) =0$. Moreover, by arguments similar to those used in the proof of
Lemma \ref{L : coeffTorus}%
\begin{equation}
\mathbb{E}\left[ \Re \left( \widetilde{a}_{\pi }\left( F\right) \right) \Im
\left( \widetilde{a}_{\pi }\left( F\right) \right) \right] =0  \label{prop0}
\end{equation}%
for every $\pi \in \hat{G}$, and also%
\begin{equation}
\widetilde{a}_{\pi }\left( F\right) =\overline{\widetilde{a}_{\pi
^{-1}}\left( F\right) }\text{ \ \ and \ \ }\Re \left( \widetilde{a}_{\pi
}\left( F\right) \right) \overset{law}{=}\Im \left( \widetilde{a}_{\pi
}\left( F\right) \right) .  \label{prop}
\end{equation}%
In general, $\Re \left( \widetilde{a}_{\pi }\left( F\right) \right) $ and $%
\Im \left( \widetilde{a}_{\pi }\left( F\right) \right) $ are not independent.

\bigskip

\textbf{Remark -- }The results of this paper extend immediately to (not
necessarily centered) functions $F\in L^{2}\left( \mathbb{R},\exp \left(
-x^{2}/2\right) dx\right) $, by considering the function $F^{\prime }=F-\int
F\left( x\right) \left( 2\pi \right) ^{-\frac{1}{2}}e^{-x^{2}/2}dx.$

\bigskip

We are interested in studying the asymptotic behavior of the coefficients $%
\widetilde{a}_{\pi }$ along some infinite sequence $\left\{ \pi _{l}:l\geq
1\right\} \subset \hat{G}$. In particular, we shall determine conditions on
the coefficients $\left\{ C_{\pi }\right\} $ in (\ref{vark}) ensuring that,
for a fixed $F$, the central limit theorem (\ref{CLT}) holds, where $%
N,N^{\prime }$ $\sim \mathcal{N}\left( 0,1/2\right) $ are independent. The
first series of results involves Hermite polynomials.

\section{Necessary and sufficient conditions for Hermite transformations
(statements)\label{S : NSher}}

We start by giving a exhaustive characterization of the CLT (\ref{CLT}),
when $F$ is an \textit{Hermite polynomial} of arbitrary order $m\geq 2$.
Recall (see e.g. \cite[p. 20]{Janss}) that the sequence $\left\{ H_{m}:m\geq
0\right\} $ of Hermite polynomials is defined through the relation

\begin{equation}
H_{m}\left( x\right) =\left( -1\right) ^{m}e^{\frac{x^{2}}{2}}\frac{d^{m}}{%
dx^{m}}\left( e^{-\frac{x^{2}}{2}}\right) \text{, \ \ }x\in \mathbb{R}\text{%
, \ }m\geq 0;  \label{Her}
\end{equation}%
it is well known that the sequence $\left\{ \left( m!\right)
^{-1/2}H_{m}:m\geq 0\right\} $ constitutes an orthonormal basis of the space
$L^{2}\left( \mathbb{R},\left( 2\pi \right) ^{-1/2}e^{-\frac{x^{2}}{2}%
}dx\right) $.

\bigskip

To state our main results, we need to introduce some further notation. For $%
\pi \in \hat{G}$ and $m\geq 1$, define the coefficient $\widehat{C}_{\pi ,m}$
as
\begin{eqnarray}
\widehat{C}_{\pi ,m} &\triangleq &\sum_{\sigma _{1}\in \hat{G}}\cdot \cdot
\cdot \sum_{\sigma _{m}\in \hat{G}}\left\{ C_{\sigma _{1}}C_{\sigma
_{2}}\cdot \cdot \cdot C_{\sigma _{m}}\right\} \mathbf{1}_{\sigma _{1}\cdot
\cdot \cdot \sigma _{m}=\pi }  \label{Clambda1} \\
&=&\sum_{\substack{ \sigma _{1},...,\sigma _{m}\in \hat{G}  \\ \sigma
_{1}\cdot \cdot \cdot \sigma _{m}=\pi }}C_{\sigma _{1}}C_{\sigma _{2}}\cdot
\cdot \cdot C_{\sigma _{m}}  \label{Clambda3} \\
&=&\sum_{\sigma _{1},...,\sigma _{m-1}\in \hat{G}}C_{\sigma _{1}}C_{\sigma
_{2}}\cdot \cdot \cdot C_{\left( \sigma _{1}\cdot \cdot \cdot \sigma
_{m-1}\right) ^{-1}\pi }.  \label{Clambda4}
\end{eqnarray}%
Note that, in (\ref{Clambda1})-(\ref{Clambda4}), $\hat{G}$ is regarded as an
Abelian group, with group operation given by (\ref{AbProd}), and that $%
\widehat{C}_{\pi ,q}=\widehat{C}_{\pi ^{-1},q}$. Moreover, $\widehat{C}_{\pi
,1}=C_{\pi }$ and, for every $m\geq 2$ and $q=1,...,m-1$,
\begin{equation}
\widehat{C}_{\pi ,m}=\sum_{\mu \in \hat{G}}\widehat{C}_{\mu ,q}\widehat{C}%
_{\pi \mu ^{-1},m-q}.  \label{Clambda5}
\end{equation}

\bigskip

In the statements of the subsequent results, we systematically adopt the
same notation and conventions pinpointed in the introduction (see formula (%
\ref{CLT})), that is: when no further specification is given, $\left\{ \pi
\right\} =\left\{ \pi _{l}:l=1,2,...\right\} $ stands for a fixed sequence
of elements of $\hat{G}$, and all limits are taken as $l\rightarrow +\infty $%
.

\bigskip

\begin{theorem}
\label{T : Main Torus}Fix $m\geq 2$, and define the random variable $%
\widetilde{a}_{\pi }\left( H_{m}\right) $ according to (\ref{FCoef}) and (%
\ref{Her}), i.e.%
\begin{equation}
\widetilde{a}_{\pi }\left( H_{m}\right) =\int_{G}H_{m}\left( T\left(
g\right) \right) \chi _{\pi }\left( g^{-1}\right) dg.  \label{FherCOEFF}
\end{equation}%
Then,
\begin{equation}
\mathbb{E}\left[ \left\vert \widetilde{a}_{\pi }\left( H_{m}\right)
\right\vert ^{2}\right] =m!\widehat{C}_{\pi ,m},  \label{Varem}
\end{equation}%
where $\widehat{C}_{\pi ,m}$ is defined as in (\ref{Clambda1}). Moreover,
the following four asymptotic conditions are equivalent:

\begin{enumerate}
\item
\begin{equation}
\frac{\widetilde{a}_{\pi }\left( H_{m}\right) }{\sqrt{m!\widehat{C}_{\pi ,m}}%
}\underset{\left\{ \pi \right\} }{\overset{law}{\rightarrow }}N+iN^{\prime },
\label{CLTher}
\end{equation}%
where $N,N^{\prime }$ $\sim \mathcal{N}\left( 0,1/2\right) $ are independent;

\item
\begin{equation}
\left[ m!\widehat{C}_{\pi ,m}\right] ^{-2}\mathbb{E}\left[ \Re \left(
\widetilde{a}_{\pi }\left( H_{m}\right) \right) ^{4}\right] \underset{%
\left\{ \pi \right\} }{\rightarrow }\frac{3}{4}\text{, \ \ and \ \ }\left[ m!%
\widehat{C}_{\pi ,m}\right] ^{-2}\mathbb{E}\left[ \Im \left( \widetilde{a}%
_{\pi }\left( H_{m}\right) \right) ^{4}\right] \underset{\left\{ \pi
\right\} }{\rightarrow }\frac{3}{4};  \label{M4clt}
\end{equation}

\item
\begin{equation}
\widehat{C}_{\pi ,m}^{-2}\sum_{\lambda \in \hat{G}}\widehat{C}_{\lambda
,q}^{2}\widehat{C}_{\pi \lambda ^{-1},m-q}^{2}\underset{\left\{ \pi \right\}
}{\rightarrow }0\text{, \ \ }\forall q=1,...,m-1\text{;}  \label{Cond2}
\end{equation}

\item
\begin{equation}
\max_{q=1,...,m-1}\frac{\sup_{\lambda \in \hat{G}}\widehat{C}_{\lambda ,q}%
\widehat{C}_{\pi \lambda ^{-1},m-q}}{\sum_{\mu \in \hat{G}}\widehat{C}_{\mu
,q}\widehat{C}_{\pi \mu ^{-1},m-q}}\underset{\left\{ \pi \right\} }{%
\rightarrow }0.  \label{Cond3}
\end{equation}
\end{enumerate}
\end{theorem}

The proof of Theorem \ref{T : Main Torus} is the object of the subsequent
sections.

\bigskip

\textbf{Remarks -- }(a) Since $H_{1}\left( x\right) =x$,
\begin{eqnarray*}
\frac{\widetilde{a}_{\pi }\left( H_{1}\right) }{\sqrt{\widehat{C}_{\pi ,1}}}
&=&\frac{\int_{G}T\left( g\right) \chi _{\pi }\left( g^{-1}\right) dg}{\sqrt{%
C_{\pi }}}=\frac{a_{\pi }}{\sqrt{C_{\pi }}}\overset{law}{=}N+iN^{\prime }, \\
N,N^{\prime } &\sim &\mathcal{N}\left( 0,1/2\right) \text{ independent,}
\end{eqnarray*}%
where we have used Lemma \ref{L : coeffTorus}, (\ref{Fcooeff}) and the fact
that $\widehat{C}_{\pi ,1}=C_{\pi }$.

(b) (\textit{An interpretation of condition }(\ref{Cond3})\textit{\ in terms
of random walks on groups}) Note $C_{\ast }\triangleq \sum_{\pi }C_{\pi }$,
and consider a sequence of independent and identically distributed $\hat{G}$%
-valued random variables $\left\{ X_{j}:j\geq 1\right\} $, such that
\begin{equation*}
\mathbb{P}\left[ X_{1}=\pi \right] =\frac{C_{\pi }}{C_{\ast }}\text{, \ \ }%
\forall \pi \in \hat{G}\text{.}
\end{equation*}%
We associate to the sequence $\left\{ X_{j}\right\} $ the $\hat{G}$-\textit{%
valued random walk} $Z=\left\{ Z_{m}:m\geq 0\right\} $, defined as $%
Z_{0}=\pi _{0}$, and $Z_{m}=X_{1}X_{2}\cdot \cdot \cdot X_{m}$ ($m\geq 1$).
Then, it is easily seen that, $\forall m\geq 2$, $\forall q=1,...,m-1$ and $%
\forall \pi \in \hat{G}$, the ratio appearing in (\ref{Cond3}) can be
rewritten as
\begin{equation*}
\frac{\sup_{\lambda \in \hat{G}}\widehat{C}_{\lambda ,q}\widehat{C}_{\pi
\lambda ^{-1},m-q}}{\sum_{\mu \in \hat{G}}\widehat{C}_{\mu ,q}\widehat{C}%
_{\pi \mu ^{-1},m-q}}=\frac{\sup_{\lambda \in \hat{G}}\mathbb{P}\left[
Z_{q}=\lambda ,Z_{m}=\pi \right] }{\mathbb{P}\left[ Z_{m}=\pi \right] }%
=\sup_{\lambda \in \hat{G}}\mathbb{P}\left[ Z_{q}=\lambda \mid Z_{m}=\pi %
\right] ,
\end{equation*}%
so that the CLT (\ref{CLTher}) holds if, and only if,
\begin{equation}
\sup_{\lambda \in \hat{G}}\mathbb{P}\left[ Z_{q}=\lambda \mid Z_{m}=\pi %
\right] \underset{\left\{ \pi \right\} }{\rightarrow }0\text{,}  \label{RWG}
\end{equation}%
for every $q=1,...,m-1.$ Condition (\ref{RWG}) can be interpreted as
follows. For every $\pi \in \hat{G}$, define a \textquotedblleft
bridge\textquotedblright\ of length $m$, from $\pi _{0}$ to $\pi $, by
conditioning $Z$ to equal $\pi $ at time $m$. Then, (\ref{RWG}) is verified
if, and only if, the probability that the bridge hits $\lambda $ at time $q$
converges to zero, uniformly on $\lambda $, as $\pi $ moves along the
sequence $\left\{ \pi \right\} $. Plainly, when (\ref{RWG}) is verified for
every $q=1,...,m-1$, one also has that
\begin{equation*}
\sup_{\lambda _{1},...,\lambda _{m-1}\in \hat{G}}\mathbb{P}\left[
Z_{1}=\lambda _{1},...,Z_{m-1}=\lambda _{m-1}\mid Z_{m}=\pi \right] \underset%
{\left\{ \pi \right\} }{\rightarrow }0\text{,}
\end{equation*}%
meaning that, asymptotically, there is no \textquotedblleft privileged
path\textquotedblright\ of length $m$\ linking $\pi _{0}$ to $\pi $.

\bigskip

Now recall that $H_{2}\left( x\right) =x^{2}-1$: by using the fact that, for
$\pi \neq \pi _{0}$, $\int_{G}\chi _{\pi }\left( g\right) dg=0$, we deduce
from Theorem \ref{T : Main Torus} the following criterion for squared
isotropic Gaussian fields on commutative groups.

\bigskip

\begin{corollary}
\label{c : Hermite 2}Let, for $\hat{G}\ni \pi \neq \pi _{0}$,%
\begin{equation*}
\widetilde{a}_{\pi }\left( H_{2}\right) =\int_{G}\left( T^{2}\left( g\right)
-1\right) \chi _{\pi }\left( g^{-1}\right) dg=\int_{G}T^{2}\left( g\right)
\chi _{\pi }\left( g^{-1}\right) dg.
\end{equation*}%
Then,
\begin{equation*}
\mathbb{E}\left[ \left\vert \widetilde{a}_{\pi }\left( H_{2}\right)
\right\vert ^{2}\right] =2\widehat{C}_{\pi ,2}=\sum_{\lambda \in \hat{G}%
}C_{\lambda }C_{\lambda ^{-1}\pi },
\end{equation*}%
and the following conditions are equivalent

\begin{enumerate}
\item
\begin{equation*}
\left( 2\widehat{C}_{\pi ,2}\right) ^{-\frac{1}{2}}\widetilde{a}_{\pi
}\left( H_{2}\right) \overset{law}{\underset{\left\{ \pi \right\} }{%
\rightarrow }}N+iN^{\prime },
\end{equation*}%
with $N,N^{\prime }$ $\sim \mathcal{N}\left( 0,1/2\right) $ independent;

\item
\begin{equation}
\frac{\sup_{\lambda \in \hat{G}}C_{\lambda }C_{\pi \lambda ^{-1}}}{\sum_{\mu
\in \hat{G}}C_{\mu }C_{\pi \mu ^{-1}}}\underset{\left\{ \pi \right\} }{%
\rightarrow }0\text{ \ }.  \label{CC}
\end{equation}
\end{enumerate}
\end{corollary}

\bigskip

Analogously, from the relation $H_{3}\left( x\right) =x^{3}-3x$ we obtain

\bigskip

\begin{corollary}
\label{c : Hermite 3}For $\pi \in \hat{G}$, let
\begin{equation*}
\widetilde{a}_{\pi }\left( H_{3}\right) =\int_{G}\left( T^{3}\left( g\right)
-3T\left( g\right) \right) \chi _{\pi }\left( g^{-1}\right) dg.
\end{equation*}%
The following conditions are equivalent

\begin{enumerate}
\item
\begin{equation*}
\left( 6\widehat{C}_{\pi ,3}\right) ^{-\frac{1}{2}}\widetilde{a}_{\pi
}\left( H_{3}\right) \overset{law}{\underset{\left\{ \pi \right\} }{%
\rightarrow }}N+iN^{\prime },
\end{equation*}%
with $N,N^{\prime }$ $\sim \mathcal{N}\left( 0,1/2\right) $ independent;

\item
\begin{equation*}
\lim_{\left\{ \pi \right\} }\frac{\sup_{\lambda \in \hat{G}}\widehat{C}%
_{\lambda ,2}C_{\pi \lambda ^{-1}}}{\sum_{\mu \in \hat{G}}\widehat{C}_{\mu
,2}C_{\pi \mu ^{-1}}}=\lim_{\left\{ \pi \right\} }\frac{\sup_{\lambda \in
\hat{G}}C_{\lambda }\widehat{C}_{\pi \lambda ^{-1},2}}{\sum_{\mu \in \hat{G}%
}C_{\mu }\widehat{C}_{\pi \mu ^{-1},2}}=0.
\end{equation*}
\end{enumerate}
\end{corollary}

\bigskip

Our strategy to prove Theorem \ref{T : Main Torus} is to represent each $%
\widetilde{a}_{\pi }\left( H_{m}\right) $ as a complex-valued functional of
a centered Gaussian measure, having the special form of a multiple Wiener-It%
\^{o} integral. To do this, we need to recall several crucial facts
concerning multiple stochastic integrals of real-valued kernels, and then to
establish some useful extensions to the case of complex-valued random
variables.

\section{Ancillary results about multiple Wiener-It\^{o} integrals \label{S
: MWII}}

In this section, we summarize some basic properties of multiple Wiener-It%
\^{o} integrals. The reader is referred e.g. to \cite[Chapter VII]{Janss}
for any explained definition or result.

\bigskip

\textit{Real kernels -- }Let $\left( A,\mathcal{A},\mu \right) $ be a finite
measure space, with $\mu $ positive, finite and non-atomic. For $d\geq 1$,
we define $L_{\mathbb{R}}^{2}\left( \mu ^{d}\right) $ and $L_{s,\mathbb{R}%
}^{2}\left( \mu ^{d}\right) $ to be the Hilbert spaces, respectively of
square-integrable, and square-integrable and symmetric real-valued functions
on $A^{d}$, with respect to the product measure $\mu ^{d}$. As usual, $L_{%
\mathbb{R}}^{2}\left( \mu ^{1}\right) =L_{s,\mathbb{R}}^{2}\left( \mu
^{1}\right) =L_{\mathbb{R}}^{2}\left( \mu \right) =L_{\mathbb{R}}^{2}\left(
A,\mathcal{A},\mu \right) $.

\bigskip

We note $\mathbf{W}=\left\{ \mathbf{W}\left( h\right) :h\in L^{2}\left( \mu
\right) \right\} $ a centered \textit{isonormal Gaussian process }over $%
L^{2}\left( \mu \right) $. This means that $\mathbf{W}$ is a centered
Gaussian family indexed by $L^{2}\left( \mu \right) $ and such that%
\begin{equation*}
\mathbb{E}\left[ \mathbf{W}\left( h^{\prime }\right) \mathbf{W}\left(
h\right) \right] =\int_{A}h^{\prime }\left( a\right) h\left( a\right) \mu
\left( da\right) \triangleq \left( h^{\prime },h\right) _{L^{2}\left( \mu
\right) }\text{,}
\end{equation*}%
for every $h,h^{\prime }\in L^{2}\left( \mu \right) $. For every $f\in L_{s,%
\mathbb{R}}^{2}\left( \mu ^{d}\right) $, we define $I_{d}\left( f\right) $
to be the the multiple Wiener-It\^{o} integral (MWII) of $f$ with respect to
$\mathbf{W}$, i.e.%
\begin{equation}
I_{d}\left( f\right) =I_{d}^{\mathbf{W}}\left( f\right) =\int_{A}\cdot \cdot
\cdot \int_{A}f\left( a_{1},...,a_{d}\right) \mathbf{W}\left( da\right)
\cdot \cdot \cdot \mathbf{W}\left( da\right) ,  \label{MWI}
\end{equation}%
where the multiple integration in (\ref{MWI}) implicitly excludes diagonals.
Recall that
\begin{equation}
\mathbb{E}\left[ I_{d}\left( f\right) I_{d^{\prime }}\left( g\right) \right]
=d!\delta _{d,d^{\prime }}\left( f,g\right) _{L_{\mathbb{R}}^{2}\left( \mu
^{d}\right) }\text{,}  \label{covMWI}
\end{equation}%
where $\delta $ is the Kronecker symbol, and therefore the application $%
f\mapsto I_{d}\left( f\right) $ defines an isomorphism between the $d$th
Wiener chaos associated to $\mathbf{W}$, and the space $L_{\mathbb{R}%
}^{2}\left( \mu ^{d}\right) $, endowed with the modified norm $\sqrt{d!}%
\left\Vert \cdot \right\Vert _{L_{\mathbb{R}}^{2}\left( \mu ^{d}\right) }$.
A fundamental relation between objects such as (\ref{MWI}) and the Hermite
polynomials introduced in (\ref{Her}) is the following: for every $h\in $ $%
L_{\mathbb{R}}^{2}\left( \mu \right) $ such that $\left\Vert h\right\Vert
_{L_{\mathbb{R}}^{2}\left( \mu \right) }=1$, and every $m\geq 1$,%
\begin{equation}
H_{m}\left( I_{1}\left( h\right) \right) =I_{m}\left( h\otimes \cdot \cdot
\cdot \otimes h\right) \text{, \ }  \label{Her>MWI}
\end{equation}%
where the tensor product inside the second integral is defined as
\begin{equation*}
h\otimes \cdot \cdot \cdot \otimes h\left( a_{1},...,a_{m}\right) =h\left(
a_{1}\right) \cdot \cdot \cdot h\left( a_{m}\right) \in L_{s,\mathbb{R}%
}^{2}\left( \mu ^{m}\right) ,
\end{equation*}%
$\forall a_{1},...,a_{m}\in A^{m}$.

\bigskip

For every $d\geq 2$, every $f\in L_{s,\mathbb{R}}^{2}\left( \mu ^{d}\right) $
and every $q=1,...,d-1$, we define the (not necessarily symmetric) \textit{%
contraction kernel} $f\otimes _{q}f\in L_{\mathbb{R}}^{2}\left( \mu
^{2\left( d-q\right) }\right) $ as%
\begin{eqnarray}
&&f\otimes _{q}f\left( x_{1},...,x_{2\left( d-q\right) }\right)
\label{contr} \\
&\triangleq &\int_{A^{q}}f\left( a_{1},...,a_{q},x_{1},...,x_{d-q}\right)
f\left( a_{1},...,a_{q},x_{d-q+1},...,x_{2\left( d-q\right) }\right) \mu
\left( da_{1}\right) \cdot \cdot \cdot \mu \left( da_{q}\right) \text{.}
\notag
\end{eqnarray}

\bigskip

The following CLT, which has been proved in \cite{NuPe} (for the Part A) and
\cite{PT} (for the Part B), concerns sequences (of vectors of) MWIIs such as
(\ref{MWI}). It is the crucial element in the proof of Theorem \ref{T : Main
Torus}.

\begin{theorem}[Nualart and Peccati, 2005; Peccati and Tudor, 2004]
\label{T : NP2005}(\textbf{A}) Fix $d\geq 2$, and let $f_{k}\in L_{s,\mathbb{%
R}}^{2}\left( \mu ^{d}\right) $, $k\geq 1$. If the variance of $I_{d}\left(
f_{k}\right) $\ converges to $1$ ($k\rightarrow +\infty $) the following
three conditions are equivalent: (i)\ $I_{d}\left( f_{k}\right) $\ converges
in law to a standard Gaussian random variable $N\left( 0,1\right) $, (ii)\ $E%
\left[ I_{d}\left( f_{k}\right) ^{4}\right] \rightarrow 3$, (iii)\ for every
$q=1,...,d-1$, the contraction kernel $f_{k}\otimes _{q}f_{k}$\ converges to
$0$ in $L_{\mathbb{R}}^{2}\left( \mu ^{2\left( d-q\right) }\right) $.

(\textbf{B}) Fix integers $p\geq 2$ and $1\leq d_{1}\leq \cdot \cdot \cdot
\leq d_{p}$. Consider a sequence of vectors
\begin{equation*}
\left( f_{k}^{\left( 1\right) },f_{k}^{\left( 2\right) },...,f_{k}^{\left(
p\right) }\right) \text{, \ \ }k\geq 1\text{,}
\end{equation*}%
such that, for each $k$, $f_{k}^{\left( j\right) }\in L_{s,\mathbb{R}%
}^{2}\left( \mu ^{d_{j}}\right) $, $j=1,...,p$, and
\begin{equation*}
\lim_{n}\mathbb{E}\left[ I_{d_{j}}\left( f_{k}^{\left( j\right) }\right)
I_{d_{i}}\left( f_{k}^{\left( i\right) }\right) \right] =\delta _{i,j}\text{,%
}
\end{equation*}%
where $\delta $ is the Kronecker symbol. Then, if $\forall $ $j=1,...,p$ the
sequence $\left\{ f_{k}^{\left( j\right) }:k\geq 1\right\} $ satisfies
either one of conditions (i)-(iii) of Part A (with $d_{j}$ substituting $d$%
), as $k\rightarrow +\infty $,
\begin{equation*}
\left( I_{d_{1}}\left( f_{k}^{\left( 1\right) }\right) ,...,I_{d_{p}}\left(
f_{k}^{\left( p\right) }\right) \right) \overset{law}{\rightarrow }\mathbf{N}%
_{p}\text{,}
\end{equation*}%
where $\mathbf{N}_{p}=\left( N_{1},...,N_{p}\right) \sim \mathcal{N}%
_{p}\left( 0,\mathbb{I}_{p}\right) $ is a $p$-dimensional vector of
independent, centered standard Gaussian random variables.
\end{theorem}

\bigskip

\textit{Complex kernels -- }For $n\geq 1$ and $d\geq 1$, $L_{\mathbb{C}%
}^{2}\left( \mu ^{d}\right) $ and $L_{s,\mathbb{C}}^{2}\left( \mu
^{d}\right) $ are the Hilbert spaces, respectively of square integrable and
square integrable and symmetric complex-valued functions with respect to the
product Lebesgue measure. For every $g\in L_{s,\mathbb{C}}^{2}\left( \mu
^{d}\right) $ with the form $g=a+ib$, where $a,b\in L_{s,\mathbb{R}%
}^{2}\left( \mu ^{d}\right) $, we set $I_{d}\left( g\right) =I_{d}\left(
a\right) +iI_{d}\left( b\right) $. Note that, by (\ref{covMWI}),
\begin{equation}
\mathbb{E}\left[ I_{d}\left( g\right) \overline{I_{d^{\prime }}\left(
f\right) }\right] =d!\delta _{d,d^{\prime }}\left( g,f\right) _{L_{\mathbb{C}%
}^{2}\left( \mu ^{d}\right) }.  \label{complcovMWI}
\end{equation}%
Also, a random variable such as $I_{d}\left( g\right) $ is real valued if,
and only if, $g\ $is real valued. For every pair and $g_{k}=a_{k}+ib_{k}\in
L_{s,\mathbb{C}}^{2}\left( \mu ^{d}\right) $, $k=1,2$, and every $%
q=1,...,d-1 $, we set
\begin{eqnarray}
&&g_{1}\otimes _{q}g_{2}\left( x_{1},...,x_{2\left( d-q\right) }\right)
\notag \\
&=&\int_{A^{q}}g_{1}\left( a_{1},...,a_{q},x_{1},...,x_{d-q}\right)
g_{2}\left( a_{1},...,a_{q},x_{d-q+1},...,x_{2\left( d-q\right) }\right) \mu
\left( da_{1}\right) \cdot \cdot \cdot \mu \left( da_{q}\right)  \notag \\
&=&a_{1}\otimes _{q}a_{2}-b_{1}\otimes _{q}b_{2}+i\left( a_{1}\otimes
_{q}b_{2}+b_{1}\otimes _{q}a_{2}\right) .  \label{complCon1}
\end{eqnarray}

The following result is an extension of Theorem \ref{T : NP2005}.

\begin{proposition}
\label{P : ConvComplx}Suppose that the sequence $g_{l}=a_{l}+ib_{l}\in L_{s,%
\mathbb{C}}^{2}\left( \mu ^{d}\right) $, $l\geq 1$, is such that
\begin{equation}
\lim_{l\rightarrow +\infty }d!\left\Vert a_{l}\right\Vert _{L_{\mathbb{R}%
}^{2}\left( \mu ^{d}\right) }^{2}=\lim_{l\rightarrow +\infty }d!\left\Vert
b_{l}\right\Vert _{L_{\mathbb{R}}^{2}\left( \mu ^{d}\right) }^{2}\rightarrow
\frac{1}{2}\text{ \ \ and \ \ }\left( a_{l},b_{l}\right) _{L_{\mathbb{R}%
}^{2}\left( \mu ^{d}\right) }=0.  \label{nor}
\end{equation}

Then, the following conditions are equivalent: as $l\rightarrow +\infty $,

\begin{enumerate}
\item $I_{d}\left( g_{l}\right) \overset{law}{\rightarrow }N+iN^{\prime }$,
where $N,N^{\prime }\sim \mathcal{N}\left( 0,1/2\right) $ are independent;

\item $g_{l}\otimes _{q}\overline{g_{l}}\rightarrow 0$ and $g_{l}\otimes
_{q}g_{l}\rightarrow 0$ in $L_{\mathbb{C}}^{2}\left( \mu ^{2\left(
d-q\right) }\right) $ for every $q=1,...,d-1$;

\item $g_{l}\otimes _{q}\overline{g_{l}}\rightarrow 0$ in $L_{\mathbb{C}%
}^{2}\left( \mu ^{2\left( d-q\right) }\right) $ for every $q=1,...,d-1$;

\item $a_{l}\otimes _{q}a_{l}\rightarrow 0$, $b_{l}\otimes
_{q}b_{l}\rightarrow 0$ and $a_{l}\otimes _{q}b_{l}\rightarrow 0$ in $L_{%
\mathbb{R}}^{2}\left( \mu ^{2\left( d-q\right) }\right) $ for every $%
q=1,...,d-1$;

\item $a_{l}\otimes _{q}a_{l}\rightarrow 0$, $b_{l}\otimes
_{q}b_{l}\rightarrow 0$ in $L_{\mathbb{R}}^{2}\left( \mu ^{2\left(
d-q\right) }\right) $ for every $q=1,...,d-1$;

\item $\mathbb{E}\left[ I_{d}\left( a_{l}\right) ^{4}\right] \rightarrow 3/4$%
, $\mathbb{E}\left[ I_{d}\left( b_{l}\right) ^{4}\right] \rightarrow 3/4$
and $\mathbb{E}\left[ I_{d}\left( a_{l}\right) ^{2}I_{d}\left( b_{l}\right)
^{2}\right] \rightarrow 1/4;$

\item $\mathbb{E}\left[ I_{d}\left( a_{l}\right) ^{4}\right] \rightarrow 3/4$%
, $\mathbb{E}\left[ I_{d}\left( b_{l}\right) ^{4}\right] \rightarrow 3/4.$
\end{enumerate}
\end{proposition}

\begin{proof}
Note first that, due to (\ref{covMWI}) and the second part of (\ref{nor}),
\begin{equation*}
\mathbb{E}\left[ I_{d}\left( b_{l}\right) I_{d}\left( a_{l}\right) \right]
=\left( a_{l},b_{l}\right) _{L_{s,\mathbb{R}}^{2}\left( \mu ^{d}\right) }=0%
\text{, \ \ }l\geq 1\text{.}
\end{equation*}%
Now, $\left( 7\rightarrow 1\right) $ holds because of (\ref{nor}) and Part B
of Theorem \ref{T : NP2005}. $\left( 5\leftrightarrow 1\rightarrow 6\right) $
is again a consequence of (\ref{nor}) and Part B of Theorem \ref{T : NP2005}
(note that (\ref{nor}) implies that all moments of the real and imaginary
parts of $I_{d}\left( a_{l}\right) $ and $I_{d}\left( b_{l}\right) $ are
uniformly bounded). $\left( 2\longleftrightarrow 4\right) $ derives from (%
\ref{complCon1}). $\left( 2\rightarrow 3\right) $, $\left( 4\rightarrow
5\right) $ and $\left( 6\rightarrow 7\right) $ and are obvious. $\left(
5\rightarrow 4\right) $ is a consequence of%
\begin{eqnarray*}
\left\Vert a_{l}\otimes _{q}b_{l}\right\Vert _{L_{\mathbb{R}}^{2}\left( \mu
^{2\left( d-q\right) }\right) }^{2}
&=&\int_{A^{d-q}}\int_{A^{d-q}}\int_{A^{q}}\int_{A^{q}}a_{l}\left( \mathbf{s}%
_{q},\mathbf{a}_{d-q}\right) b_{l}\left( \mathbf{s}_{q},\mathbf{b}%
_{d-q}\right) a_{l}\left( \mathbf{t}_{q},\mathbf{a}_{d-q}\right) \\
&&b_{l}\left( \mathbf{t}_{q},\mathbf{b}_{d-q}\right) \mu ^{d-q}\left( d%
\mathbf{a}_{d-q}\right) \mu ^{d-q}\left( d\mathbf{b}_{d-q}\right) \mu
^{q}\left( d\mathbf{s}_{q}\right) \mu ^{q}\left( d\mathbf{t}_{q}\right) \\
&=&\left( \left( a_{l}\otimes _{d-q}a_{l}\right) ,\left( b_{l}\otimes
_{d-q}b_{l}\right) \right) _{L_{\mathbb{R}}^{2}\left( \mu ^{2q}\right) },
\end{eqnarray*}%
where $\mathbf{s}_{q}$ stands for a vector of the type
\begin{equation*}
\left( s_{1},...,s_{q}\right) \text{, \ \ \ with }s_{j}\in A\text{, \ \ }%
j=1,...,q\text{,}
\end{equation*}%
and $\mu ^{q}\left( d\mathbf{s}_{q}\right) =\mu \left( ds_{1}\right) \cdot
\cdot \cdot \mu \left( ds_{q}\right) $ (similar conventions apply to $%
\mathbf{a}_{d-q}$, $\mathbf{b}_{d-q}$ and $\mathbf{t}_{q}$). We are left
with the implication $\left( 3\rightarrow 2\right) $, which is a consequence
of the relation%
\begin{equation}
\left\Vert g_{l}\otimes _{q}\overline{g_{l}}\right\Vert _{L_{\mathbb{R}%
}^{2}\left( \mu ^{2\left( d-q\right) }\right) }^{2}\geq \left\Vert
g_{l}\otimes _{q}g_{l}\right\Vert _{L_{\mathbb{R}}^{2}\left( \mu ^{2\left(
d-q\right) }\right) }^{2}\text{, \ \ }\forall l\geq 1\text{.}  \label{o}
\end{equation}

To prove (\ref{o}), just write
\begin{eqnarray*}
\left\Vert g_{l}\otimes _{q}\overline{g_{l}}\right\Vert _{L_{\mathbb{R}%
}^{2}\left( \mu ^{2\left( d-q\right) }\right) }^{2} &=&\left\Vert
a_{l}\otimes _{q}a_{l}\right\Vert _{L_{\mathbb{R}}^{2}\left( \mu ^{2\left(
d-q\right) }\right) }^{2}+\left\Vert b_{l}\otimes _{q}b_{l}\right\Vert _{L_{%
\mathbb{R}}^{2}\left( \mu ^{2\left( d-q\right) }\right) }^{2} \\
&&+2\left( a_{l}\otimes _{q}a_{l},b_{l}\otimes _{q}b_{l}\right) _{L_{\mathbb{%
R}}^{2}\left( \mu ^{2\left( d-q\right) }\right) }+2\left\Vert a_{l}\otimes
_{q}b_{l}\right\Vert _{L_{\mathbb{R}}^{2}\left( \mu ^{2\left( d-q\right)
}\right) }^{2} \\
&&-2\left( a_{l}\otimes _{q}b_{l},b_{l}\otimes _{q}a_{l}\right) _{L_{\mathbb{%
R}}^{2}\left( \mu ^{2\left( d-q\right) }\right) }
\end{eqnarray*}%
and%
\begin{eqnarray*}
\left\Vert g_{l}\otimes _{q}g_{l}\right\Vert _{L_{\mathbb{R}}^{2}\left( \mu
^{2\left( d-q\right) }\right) }^{2} &=&\left\Vert a_{l}\otimes
_{q}a_{l}\right\Vert _{L_{\mathbb{R}}^{2}\left( \mu ^{2\left( d-q\right)
}\right) }^{2}+\left\Vert b_{l}\otimes _{q}b_{l}\right\Vert _{L_{\mathbb{R}%
}^{2}\left( \mu ^{2\left( d-q\right) }\right) }^{2} \\
&&-2\left( a_{l}\otimes _{q}a_{l},b_{l}\otimes _{q}b_{l}\right) _{L_{\mathbb{%
R}}^{2}\left( \mu ^{2\left( d-q\right) }\right) }+2\left\Vert a_{l}\otimes
_{q}b_{l}\right\Vert _{L_{\mathbb{R}}^{2}\left( \mu ^{2\left( d-q\right)
}\right) }^{2} \\
&&+2\left( a_{l}\otimes _{q}b_{l},b_{l}\otimes _{q}a_{l}\right) _{L_{\mathbb{%
R}}^{2}\left( \mu ^{2\left( d-q\right) }\right) },
\end{eqnarray*}%
and finally
\begin{eqnarray*}
&&2\left( a_{l}\otimes _{q}a_{l},b_{l}\otimes _{q}b_{l}\right) _{L_{\mathbb{R%
}}^{2}\left( \mu ^{2\left( d-q\right) }\right) }-2\left( a_{l}\otimes
_{q}b_{l},b_{l}\otimes _{q}a_{l}\right) _{L_{\mathbb{R}}^{2}\left( \mu
^{2\left( d-q\right) }\right) } \\
&=&2\left( a_{l}\otimes _{d-q}b_{l},a_{l}\otimes _{d-q}b_{l}\right) _{L_{%
\mathbb{R}}^{2}\left( \mu ^{2q}\right) }-2\left( a_{l}\otimes
_{d-q}b_{l},b_{l}\otimes _{d-q}a_{l}\right) _{L_{\mathbb{R}}^{2}\left( \mu
^{2q}\right) } \\
&=&\left\Vert a_{l}\otimes _{d-q}b_{l}\right\Vert _{L_{\mathbb{R}}^{2}\left(
\mu ^{2q}\right) }^{2}+\left\Vert b_{l}\otimes _{d-q}a_{l}\right\Vert _{L_{%
\mathbb{R}}^{2}\left( \mu ^{2q}\right) }^{2}-2\left( a_{l}\otimes
_{d-q}b_{l},b_{l}\otimes _{d-q}a_{l}\right) _{L_{\mathbb{R}}^{2}\left( \mu
^{2q}\right) } \\
&=&\left\Vert a_{l}\otimes _{d-q}b_{l}-b_{l}\otimes _{d-q}a_{l}\right\Vert
_{L_{\mathbb{R}}^{2}\left( \mu ^{2q}\right) }^{2}\geq 0.
\end{eqnarray*}
\end{proof}

\section{Proof of Theorem \protect\ref{T : Main Torus} \label{S : Proof}}

Let $\left\{ C_{\pi }:\pi \in \hat{G}\right\} $ be defined as in (\ref{vark}%
). We start by considering a collection of complex-valued and square
integrable functions $\left\{ f_{\pi }:\pi \in \hat{G}\right\} \subset L_{%
\mathbb{C}}^{2}\left( \mu \right) $, with the following properties: (i) $\Im
\left( f_{\pi _{0}}\right) =0$, (ii) $f_{\pi }=\overline{f}_{\pi ^{-1}}$,
(iii) $\int_{A}f_{\pi }\left( a\right) f_{\pi }\left( a\right) \mu \left(
da\right) =0$, $\forall \pi \neq \pi _{0}$, (iv) both $\Re \left( f_{\pi
}\right) $ and $\Im \left( f_{\pi }\right) $ are orthogonal (in $L_{\mathbb{R%
}}^{2}\left( \mu \right) $) to $\Re \left( f_{\sigma }\right) $ and $\Im
\left( f_{\sigma }\right) $ for every $\sigma \notin \left\{ \pi ,\pi
^{-1}\right\} $, (v) $\int_{A}\left\vert f_{\pi }\left( a\right) \right\vert
^{2}\mu \left( da\right) =C_{\pi }$. Note that
\begin{equation*}
\int_{A}f_{\pi }\left( a\right) f_{\pi }\left( a\right) \mu \left( da\right)
=\int_{A}\left( \Re \left( f_{\pi }\left( a\right) \right) ^{2}-\Im \left(
f_{\pi }\left( a\right) \right) ^{2}\right) \mu \left( da\right)
+2i\int_{A}\Re \left( f_{\pi }\left( a\right) \right) \Im \left( f_{\pi
}\left( a\right) \right) \mu \left( da\right) ,
\end{equation*}%
and therefore (iii) holds if, and only if,
\begin{equation*}
C_{\pi }=\int_{A}\left\vert f_{\pi }\left( a\right) \right\vert ^{2}\mu
\left( da\right) =2\int_{A}\Re \left( f_{\pi }\left( a\right) \right)
^{2}\mu \left( da\right) =2\int_{A}\Im \left( f_{\pi }\left( a\right)
\right) ^{2}\mu \left( da\right) ,\text{ \ \ }\forall \pi \neq \pi _{0}\text{%
,}
\end{equation*}%
and $\int_{A}\Re \left( f_{\pi }\left( a\right) \right) \Im \left( f_{\pi
}\left( a\right) \right) \mu \left( da\right) =0$ for every $\pi .$

\bigskip

The class $\left\{ f_{\pi }:\pi \in \hat{G}\right\} $ can be constructed as
follows. Let $\left\{ ...,\pi _{-1},\pi _{0},\pi _{1},\pi _{2},...\right\} $
be any two-sided enumeration of $\hat{G}$, such that $\pi _{0}$ is the
trivial representation as before, and $\pi _{j}=\pi _{-j}^{-1}$ for every $%
j=1,2,...$. Then, consider an orthonormal basis $\left\{
e_{k}:k=...,-1,0,1,2,...\right\} $ of $L_{\mathbb{R}}^{2}\left( \mu \right) $%
, and set $f_{\pi _{0}}=e_{0}$ and, for $j\geq 1$,%
\begin{equation*}
f_{\pi _{j}}=\sqrt{\frac{C_{\pi _{j}}}{2}}\times \left( e_{j}+ie_{-j}\right)
\text{ \ and \ }f_{\pi _{-j}}=\sqrt{\frac{C_{\pi _{j}}}{2}}\times \left(
e_{j}-ie_{-j}\right)
\end{equation*}%
(with this notation, one has plainly that $C_{\pi _{j}}=C_{\pi _{-j}}$).

\bigskip

The next Lemma is easily verified.

\begin{lemma}
\label{L : Identity}The following identity in law holds
\begin{equation}
\left\{ I_{1}\left( f_{\pi }\right) :\pi \in \hat{G}\right\} \overset{law}{=}%
\left\{ a_{\pi }:\pi \in \hat{G}\right\} \text{,}  \label{IDL}
\end{equation}%
where the coefficients $a_{\pi }$ are given by (\ref{Fcooeff}), and therefore%
\begin{equation}
T\left( g\right) \overset{law}{=}\sum_{\pi \in \hat{G}}I_{1}\left( f_{\pi
}\right) \chi _{\pi }\left( g\right) \text{, \ \ }g\in G\text{,}  \label{IL}
\end{equation}%
where the identity in law is in the sense of stochastic processes. As a
consequence, for every $F\in L^{2}\left( \mathbb{R},\exp \left(
-x^{2}/2\right) dx\right) $ and $g\in G$%
\begin{equation}
\widetilde{a}_{\pi }\left( F\right) \overset{law}{=}\int_{G}F\left[
\sum_{\pi \in \hat{G}}I_{1}\left( f_{\pi }\right) \chi _{\pi }\left(
g\right) \right] \chi _{\pi }\left( g^{-1}\right) dg\text{,}  \label{I}
\end{equation}%
where $\widetilde{a}_{\pi }\left( F\right) $ is defined as in (\ref{FCoef}).
\end{lemma}

\bigskip

Since, for any $\pi \in G$,%
\begin{equation*}
\sum_{\pi \in \hat{G}}C_{\pi }=\sum_{\pi \in \hat{G}}\left\Vert f_{\pi
}\right\Vert _{L_{\mathbb{C}}^{2}\left( \mu \right) }^{2}=\sum_{\pi \in \hat{%
G}}\left\Vert f_{\pi }\chi _{\pi }\left( g\right) \right\Vert _{L_{\mathbb{C}%
}^{2}\left( \mu \right) }^{2}<+\infty ,
\end{equation*}%
for every \textit{fixed} $g\in G$ and any sequence of finite subsets $\hat{G}%
_{N}\subset \hat{G}$ such that $\hat{G}_{N}\uparrow \hat{G}$, the sequence%
\begin{equation*}
\sum_{g\in \hat{G}_{N}}f_{\pi }\left( \cdot \right) \chi _{\pi }\left(
g\right) \in L_{\mathbb{C}}^{2}\left( \mu \right) \text{, \ \ }N\geq 1\text{,%
}
\end{equation*}%
converges (as $N\rightarrow +\infty $) in $L_{\mathbb{C}}^{2}\left( \mu
\right) $ to a certain function
\begin{equation}
h^{g}\left( \cdot \right) \triangleq \sum_{\pi \in \hat{G}}f_{\pi }\left(
\cdot \right) \chi _{\pi }\left( g\right) \in L_{\mathbb{C}}^{2}\left( \mu
\right)  \label{pd0}
\end{equation}%
(we stress that in (\ref{pd0}) $g$ is a fixed parameter). Note that the
properties of the $f_{\pi }$'s imply that $h^{g}$ is real-valued, and also
that the mapping $\left( x,g\right) \mapsto h^{g}\left( x\right) $ is
jointly measurable. By using the linearity of MWIIs, we deduce from (\ref{IL}%
) that, as stochastic processes,%
\begin{equation}
T\left( g\right) \overset{law}{=}I_{1}\left( h^{g}\right) \text{, \ \ }g\in G%
\text{,}  \label{pd1}
\end{equation}%
and therefore (\ref{I}) implies that for every $\pi $,
\begin{equation}
\widetilde{a}_{\pi }\left( F\right) \overset{law}{=}\int_{G}F\left[
I_{1}\left( h^{g}\right) \right] \chi _{\pi }\left( g\right) dg.  \label{pd2}
\end{equation}

\bigskip

Now fix $m\geq 2$, and consider the $m$th Hermite polynomial $H_{m}$. Since
\begin{equation*}
1=\mathbb{E}\left[ T\left( g\right) ^{2}\right] =\mathbb{E}\left[
I_{1}\left( h^{g}\right) ^{2}\right] =\left\Vert h^{g}\right\Vert _{L_{%
\mathbb{R}}^{2}\left( \mu \right) }^{2},
\end{equation*}%
we deduce from (\ref{Her>MWI}) that, for every $g\in G$, $H_{m}\left[
I_{1}\left( h^{g}\right) \right] =I_{m}\left( h^{g}\otimes \cdot \cdot \cdot
\otimes h^{g}\right) $. Thus, by using (\ref{pd2}) in the case $F=H_{m}$ and
by interchanging deterministic and stochastic integration,
\begin{eqnarray}
&&\widetilde{a}_{\pi }\left( H_{m}\right) \overset{law}{=}\int_{G}H_{m}\left[
I_{1}\left( h^{g}\right) \right] \chi _{\pi }\left( g^{-1}\right) dg
\label{I1} \\
&=&\int_{G}I_{m}\left( h^{g}\otimes \cdot \cdot \cdot \otimes h^{g}\right)
\chi _{\pi }\left( g^{-1}\right) dg  \notag \\
&=&I_{m}\left( \int_{G}\left\{ h^{g}\otimes \cdot \cdot \cdot \otimes
h^{g}\right\} \chi _{\pi }\left( g^{-1}\right) dg\right)  \label{SFT} \\
&=&I_{m}\left( \widetilde{h}_{m,\pi }\right) ,  \notag
\end{eqnarray}%
where
\begin{equation}
\widetilde{h}_{m,\pi }\triangleq \int_{G}\left\{ h^{g}\otimes \cdot \cdot
\cdot \otimes h^{g}\right\} \chi _{\pi }\left( g^{-1}\right) dg\in L_{s,%
\mathbb{C}}^{2}\left( \mu ^{m}\right) .  \label{kerdef}
\end{equation}

\bigskip

\textbf{Remark -- }Since the Haar measure $dg$ has finite mass, the
\textquotedblleft stochastic Fubini theorem\textquotedblright\ applied in (%
\ref{SFT}) can be justified by standard arguments. See for instance \cite[%
Lemma 13]{Pec2001}.

\bigskip

The function $\widetilde{h}_{m,\pi }$ can be made explicit by means of (\ref%
{pd0}), i.e.
\begin{eqnarray}
&&\widetilde{h}_{m,\pi }\left( x_{1},...,x_{m}\right)  \notag \\
&=&\int_{G}\left\{ \sum_{\sigma _{1}\in \hat{G}}f_{\sigma _{1}}\left(
x_{1}\right) \chi _{\sigma _{1}}\left( g\right) \times \cdot \cdot \cdot
\times \sum_{\sigma _{m}\in \hat{G}}f_{\sigma _{m}}\left( x_{m}\right) \chi
_{\pi _{m}}\left( g\right) \right\} \chi _{\pi }\left( g^{-1}\right) dg
\label{tildkernel} \\
&=&\sum_{\sigma _{1}\in \hat{G}}\sum_{\sigma _{2}\in \hat{G}}\cdot \cdot
\cdot \sum_{\sigma _{m-1}\in \hat{G}}f_{\sigma _{1}}\left( x_{1}\right)
f_{\sigma _{2}}\left( x_{2}\right) \times \cdot \cdot \cdot \times f_{\pi
\left( \sigma _{1}\cdot \cdot \cdot \sigma _{m-1}\right) ^{-1}}\left(
x_{m}\right) ,\text{ \ \ }\left( x_{1},...,x_{m}\right) \in A^{m}\text{,}
\notag
\end{eqnarray}%
where we used (\ref{AbCar}) and the orthogonality between characters of
non-equivalent representations. By using (\ref{I1}) and (\ref{complcovMWI}),
\begin{equation*}
\mathbb{E}\left[ \left\vert \widetilde{a}_{\pi }\left( H_{m}\right)
\right\vert ^{2}\right] =m!\left\Vert \widetilde{h}_{m,\pi }\right\Vert _{L_{%
\mathbb{C}}^{2}\left( \mu ^{m}\right) }^{2}=m!\sum_{\substack{ \sigma
_{1},...,\sigma _{m}\in \hat{G}  \\ \sigma _{1}\cdot \cdot \cdot \sigma
_{m}=\pi }}C_{\sigma _{1}}C_{\sigma _{2}}\cdot \cdot \cdot C_{\sigma _{m}}=m!%
\widehat{C}_{\pi ,m}\text{,}
\end{equation*}%
thus proving (\ref{Varem}). Now define
\begin{equation}
\widetilde{\widetilde{a}}_{\pi }\left( H_{m}\right) \triangleq \frac{%
\widetilde{a}_{\pi }\left( H_{m}\right) }{\mathbb{E}\left[ \left\vert
\widetilde{a}_{\pi }\left( H_{m}\right) \right\vert ^{2}\right] ^{\frac{1}{2}%
}}\overset{law}{=}I_{m}\left( \widetilde{\widetilde{h}}_{m,\pi }\right)
\text{,}  \label{deff}
\end{equation}%
where
\begin{equation}
\widetilde{\widetilde{h}}_{m,\pi }\triangleq \mathbb{E}\left[ \left\vert
\widetilde{a}_{\pi }\left( H_{m}\right) \right\vert ^{2}\right] ^{-\frac{1}{2%
}}\widetilde{h}_{m,\pi }=\left( m!\widehat{C}_{\pi ,m}\right) ^{-1/2}%
\widetilde{h}_{m,\pi }.  \label{deff2}
\end{equation}%
Since (\ref{prop0}) and (\ref{prop}) hold (with $F=H_{m}$), it is clear
that, for $\pi \in \hat{G}$,
\begin{equation*}
m!\left( \Re \left( \widetilde{\widetilde{h}}_{m,\pi }\right) ,\Im \left(
\widetilde{\widetilde{h}}_{m,\pi }\right) \right) _{L_{\mathbb{R}}^{2}\left(
\mu ^{m}\right) }=\mathbb{E}\left[ \Re \left( \widetilde{\widetilde{a}}_{\pi
}\left( H_{m}\right) \right) \Im \left( \widetilde{\widetilde{a}}_{\pi
}\left( H_{m}\right) \right) \right] =0
\end{equation*}%
and also $m!\left\Vert \widetilde{\widetilde{h}}_{m,\pi }\right\Vert _{L_{%
\mathbb{C}}^{2}\left( \mu ^{m}\right) }^{2}=\mathbb{E}\left[ \left\vert
\left( \widetilde{\widetilde{a}}_{\pi }\left( H_{m}\right) \right)
\right\vert ^{2}\right] =1$, so that%
\begin{eqnarray*}
\mathbb{E}\left[ \Re \left( \widetilde{\widetilde{a}}_{\pi }\left(
H_{m}\right) \right) ^{2}\right] &=&\mathbb{E}\left[ \Im \left( \widetilde{%
\widetilde{a}}_{\pi }\left( H_{m}\right) \right) ^{2}\right] =m!\left\Vert
\Re \left( \widetilde{\widetilde{h}}_{m,\pi }\right) \right\Vert _{L_{%
\mathbb{C}}^{2}\left( \mu ^{m}\right) }^{2} \\
&=&m!\left\Vert \Im \left( \widetilde{\widetilde{h}}_{m,\pi }\right)
\right\Vert _{L_{\mathbb{C}}^{2}\left( \mu ^{m}\right) }^{2}=\frac{1}{2}.
\end{eqnarray*}

It follows that all the assumptions of Proposition \ref{P : ConvComplx} are
satisfied, with $d=m$, $g_{l}=\widetilde{\widetilde{h}}_{m,\pi _{l}}$, and
therefore $a_{l}=\Re \left( \widetilde{\widetilde{h}}_{m,\pi _{l}}\right) $
and $b_{n,l}=\Im \left( \widetilde{\widetilde{h}}_{m,\pi _{l}}\right) $
(recall that, in the statement of Theorem \ref{T : Main Torus}, $\left\{ \pi
\right\} $ stands for a sequence of the form $\left\{ \pi _{l}:l\geq
1\right\} $). As a consequence, in view of (\ref{deff}), we deduce from the
implications $\left( 1\leftrightarrow 7\right) $ in Proposition \ref{P :
ConvComplx} that the convergence in law (\ref{CLTher}) holds if, and only
if, (\ref{M4clt}) is verified. We have therefore proved that Conditions 1
and 2 in Theorem \ref{T : Main Torus} are equivalent.

\bigskip

To conclude the proof, we start by observing that, thanks e.g. to the
implications $\left( 7\longleftrightarrow 1\longleftrightarrow 3\right) $ in
Proposition \ref{P : ConvComplx}, either one of conditions (\ref{CLTher})
and (\ref{M4clt}) is equivalent to the following:
\begin{equation}
\widetilde{\widetilde{h}}_{m,\pi }\otimes _{q}\overline{\left( \widetilde{%
\widetilde{h}}_{m,\pi }\right) }\underset{\left\{ \pi \right\} }{\rightarrow
}0\text{, \ in \ }L_{\mathbb{C}}^{2}\left( \mu ^{2\left( m-q\right) }\right)
,\text{ \ }\forall q\in \left\{ 1,...,m-1\right\} \text{.}  \label{cvcntrs}
\end{equation}%
It follows that the equivalence of Conditions 1, 2 and 3 in Theorem \ref{T :
Main Torus} is established, once it is shown that (\ref{cvcntrs}) is true
if, and only if, condition (\ref{Cond2}) is verified for every $q=1,...,m-1$%
. Start with $q=m-1$. Indeed,
\begin{eqnarray*}
&&\widetilde{\widetilde{h}}_{m,\pi }\otimes _{m-1}\overline{\left(
\widetilde{\widetilde{h}}_{m,\pi }\right) }\left( x_{1},x_{2}\right) \\
&=&\left( m!\widehat{C}_{\pi ,m}\right) ^{-1}\int_{A}\cdot \cdot \cdot
\int_{A}\widetilde{h}_{m,\pi }\left( a_{1},...,a_{m-1},x_{1}\right)
\overline{\left( \widetilde{h}_{m,\pi }\right) }\left(
a_{1},...,a_{m-1},x_{2}\right) \mu \left( da_{1}\right) ...\mu \left(
da_{m}\right) \\
&=&\left( m!\widehat{C}_{\pi ,m}\right) ^{-1}\sum_{\pi _{1}\in \hat{G}%
}\sum_{\pi _{2}\in \hat{G}}\cdot \cdot \cdot \sum_{\pi _{m-1}\in \hat{G}%
}C_{\pi _{1}}C_{\pi _{2}}\cdot \cdot \cdot C_{\pi _{m-1}}\times f_{\pi
\left( \pi _{1}\cdot \cdot \cdot \pi _{m-1}\right) ^{-1}}\left( x_{1}\right)
f_{\pi ^{-1}\left( \pi _{1}\cdot \cdot \cdot \pi _{m-1}\right) }\left(
x_{2}\right) \\
&=&\left( m!\widehat{C}_{\pi ,m}\right) ^{-1}\sum_{\lambda \in \hat{G}}\sum
_{\substack{ \pi _{1},...,\pi _{m-1}\in \hat{G}  \\ \pi _{1}\cdot \cdot
\cdot \pi _{m-1}=\lambda }}C_{\pi _{1}}C_{\pi _{2}}\cdot \cdot \cdot C_{\pi
_{m-1}}\times f_{\pi \lambda ^{-1}}\left( x_{1}\right) f_{\pi ^{-1}\lambda
}\left( x_{2}\right) \\
&=&\left( m!\widehat{C}_{\pi ,m}\right) ^{-1}\sum_{\lambda \in \hat{G}}%
\widehat{C}_{\lambda ,m-1}f_{\pi \lambda ^{-1}}\left( x_{1}\right) f_{\pi
^{-1}\lambda }\left( x_{2}\right) =\left( m!\widehat{C}_{\pi ,m}\right)
^{-1}\sum_{\lambda \in \hat{G}}\widehat{C}_{\lambda ,m-1}f_{\pi \lambda
^{-1}}\left( x_{1}\right) \overline{f_{\pi \lambda ^{-1}}\left( x_{2}\right)
},
\end{eqnarray*}%
yielding
\begin{equation*}
\left\Vert \widetilde{\widetilde{h}}_{m,\pi }\otimes _{m-1}\overline{\left(
\widetilde{\widetilde{h}}_{m,\pi }\right) }\right\Vert _{L_{\mathbb{C}%
}^{2}\left( \mu ^{2}\right) }^{2}=\left( m!\widehat{C}_{\pi ,m}\right)
^{-2}\sum_{\lambda \in \hat{G}}\widehat{C}_{\lambda ,m-1}^{2}C_{\pi \lambda
^{-1}}C_{\pi ^{-1}\lambda }=\left( m!\widehat{C}_{\pi ,m}\right)
^{-2}\sum_{\lambda \in \hat{G}}\widehat{C}_{\pi ,m-1}^{2}C_{\pi \lambda
^{-1}}^{2}\text{,}
\end{equation*}%
thus proving that (\ref{Cond2}) holds for $q=m-1$ if, and only if, $%
\widetilde{\widetilde{h}}_{m,\pi }\otimes _{m-1}\overline{\left( \widetilde{%
\widetilde{h}}_{m,\pi }\right) }\underset{\left\{ \pi \right\} }{\rightarrow
}0$. Now suppose $m\geq 3$, and fix $q=1,...,m-2$. In this case,
\begin{eqnarray*}
&&\widetilde{\widetilde{h}}_{m,\pi }\otimes _{q}\overline{\left( \widetilde{%
\widetilde{h}}_{m,\pi }\right) }\left( x_{1},...,x_{2\left( m-q\right)
}\right) \\
&=&\left( m!\widehat{C}_{\pi ,m}\right) ^{-1}\int_{A}\cdot \cdot \cdot
\int_{A}\widetilde{h}_{m,\pi }\left(
a_{1},...,a_{q},x_{1},...,x_{m-q}\right) \times \\
&&\text{ \ \ \ \ \ \ \ \ \ \ \ \ \ \ }\times \overline{\left( \widetilde{h}%
_{m,\pi }\right) }\left( a_{1},...,a_{q},x_{m-q+1},...,x_{2\left( m-q\right)
}\right) \mu \left( da_{1}\right) \cdot \cdot \cdot \mu \left( da_{q}\right)
\\
&=&\left( m!\widehat{C}_{\pi ,m}\right) ^{-1}\sum_{\pi _{1},...,\pi _{q}\in
\hat{G}}C_{\pi _{1}}\cdot \cdot \cdot C_{\pi _{q}}\times \\
&&\text{ \ \ \ \ \ \ \ \ \ \ \ \ \ \ \ \ }\times \sum_{\rho _{1},...,\rho
_{m-1-q}}\sum_{\sigma _{1},...,\sigma _{m-1-q}}\prod_{r=1}^{m-q-1}f_{\rho
_{r}}\left( x_{r}\right) f_{\sigma _{r}^{-1}}\left( x_{m-q+r}\right) \\
&&\text{ \ \ \ \ \ \ \ \ \ \ \ \ \ \ \ \ \ \ \ \ \ \ \ \ \ \ \ }\times
f_{\pi \left( \pi _{1}\cdot \cdot \cdot \pi _{q}\right) ^{-1}\left( \rho
_{1}\cdot \cdot \cdot \rho _{m-1-q}\right) ^{-1}}\left( x_{m-q}\right)
f_{\pi ^{-1}\left( \pi _{1}\cdot \cdot \cdot \pi _{q}\right) \left( \sigma
_{1}\cdot \cdot \cdot \sigma _{m-1-q}\right) }\left( x_{2\left( m-q\right)
}\right) \\
&=&\left( m!\widehat{C}_{\pi ,m}\right) ^{-1}\sum_{\rho _{1},...,\rho
_{m-1-q}}\sum_{\sigma _{1},...,\sigma _{m-1-q}}\prod_{r=1}^{m-q-1}f_{\rho
_{r}}\left( x_{r}\right) f_{\sigma _{r}^{-1}}\left( x_{m-q+r}\right) \\
&&\text{ \ \ \ \ \ \ \ \ \ \ \ \ \ \ \ \ \ \ \ \ \ \ \ }\times \sum_{\lambda
\in \hat{G}}\widehat{C}_{\lambda ,q}f_{\pi \lambda ^{-1}\left( \rho
_{1}\cdot \cdot \cdot \rho _{m-1-q}\right) ^{-1}}\left( x_{m-q}\right)
f_{\pi ^{-1}\lambda \left( \sigma _{1}\cdot \cdot \cdot \sigma
_{m-1-q}\right) }\left( x_{2\left( m-q\right) }\right) ,
\end{eqnarray*}%
and some calculations yield%
\begin{equation}
\left\Vert \widetilde{\widetilde{h}}_{m,\pi }\otimes _{q}\overline{\left(
\widetilde{\widetilde{h}}_{m,\pi }\right) }\right\Vert _{L_{\mathbb{C}%
}^{2}\left( \mu ^{2\left( m-q\right) }\right) }^{2}=\left( m!\widehat{C}%
_{\pi ,m}\right) ^{-2}\sum_{\lambda \in \hat{G}}\widehat{C}_{\lambda ,q}^{2}%
\widehat{C}_{\pi \lambda ^{-1},m-q}^{2}\text{.}  \label{ccc}
\end{equation}

Relation (\ref{ccc}) shows in particular that, for $q=1,...,m-2$, $%
\widetilde{\widetilde{h}}_{m,\pi }\otimes _{q}\overline{\left( \widetilde{%
\widetilde{h}}_{m,\pi }\right) }\underset{\left\{ \pi \right\} }{\rightarrow
}0$ if, and only if, (\ref{Cond2}) is verified. To see that Conditions 3 and
4 in the statement of Theorem \ref{T : Main Torus} are equivalent, use (\ref%
{Clambda5}) to write
\begin{eqnarray*}
\widehat{C}_{\pi ,m}^{-2}\sum_{\lambda \in \hat{G}}\widehat{C}_{\lambda
,q}^{2}\widehat{C}_{\pi \lambda ^{-1},m-q}^{2} &=&\sum_{\lambda \in \hat{G}%
}\left\{ \frac{\widehat{C}_{\lambda ,q}\widehat{C}_{\pi \lambda ^{-1},m-q}}{%
\sum_{\mu \in \hat{G}}\widehat{C}_{\mu ,q}\widehat{C}_{\pi \mu ^{-1},m-q}}%
\right\} ^{2} \\
&\leq &\sup_{\lambda \in \hat{G}}\left\{ \frac{\widehat{C}_{\lambda ,q}%
\widehat{C}_{\pi \lambda ^{-1},m-q}}{\sum_{\mu \in \hat{G}}\widehat{C}_{\mu
,q}\widehat{C}_{\pi \mu ^{-1},m-q}}\right\} \sum_{\lambda \in \hat{G}%
}\left\{ \frac{\widehat{C}_{\lambda ,q}\widehat{C}_{\pi \lambda ^{-1},m-q}}{%
\sum_{\mu \in \hat{G}}\widehat{C}_{\mu ,q}\widehat{C}_{\pi \mu ^{-1},m-q}}%
\right\} \\
&=&\sup_{\lambda \in \hat{G}}\left\{ \frac{\widehat{C}_{\lambda ,q}\widehat{C%
}_{\pi \lambda ^{-1},m-q}}{\sum_{\mu \in \hat{G}}\widehat{C}_{\mu ,q}%
\widehat{C}_{\pi \mu ^{-1},m-q}}\right\} ,
\end{eqnarray*}%
and also
\begin{eqnarray*}
\max_{q=1,...,m-1}\widehat{C}_{\pi ,m}^{-2}\sum_{\lambda \in \hat{G}}%
\widehat{C}_{\lambda ,q}^{2}\widehat{C}_{\pi \lambda ^{-1},m-q}^{2}
&=&\max_{q=1,...,m-1}\sum_{\lambda \in \hat{G}}\left\{ \frac{\widehat{C}%
_{\lambda ,q}\widehat{C}_{\pi \lambda ^{-1},m-q}}{\sum_{\mu \in \hat{G}}%
\widehat{C}_{\mu ,q}\widehat{C}_{\pi \mu ^{-1},m-q}}\right\} ^{2} \\
&\geq &\max_{q=1,...,m-1}\sup_{\lambda \in \hat{G}}\left\{ \frac{\widehat{C}%
_{\lambda ,q}\widehat{C}_{\pi \lambda ^{-1},m-q}}{\sum_{\mu \in \hat{G}}%
\widehat{C}_{\mu ,q}\widehat{C}_{\pi \mu ^{-1},m-q}}\right\} ^{2}.
\end{eqnarray*}%
This concludes the proof of Theorem \ref{T : Main Torus}. $\ \blacksquare $

\bigskip

In Section \ref{S : genCLT}, we will establish a CLT of the type (\ref{CLT})
for functions $F\in L_{0}^{2}\left( \mathbb{R},e^{-x^{2}/2}dx\right) $ that
are not necessarily Hermite polynomials. As a first step, in the next
section we prove a result concerning the joint convergence of vectors of
coefficients of the type $\widetilde{a}_{\pi }\left( H_{m}\right) $.

\section{Joint convergence of the $\widetilde{\widetilde{a}}_{\protect\pi %
}\left( H_{m}\right) $}

Fix integers $p\geq 2$ and $2\leq m_{1}<\cdot \cdot \cdot <m_{p}$, and
define, for $\pi \in \hat{G}$, the vectors
\begin{equation*}
\left( \widetilde{a}_{\pi }\left( H_{m_{1}}\right) ,...,\widetilde{a}_{\pi
}\left( H_{m_{p}}\right) \right) \text{ \ \ and \ \ }\left( \widetilde{%
\widetilde{a}}_{\pi }\left( H_{m_{1}}\right) ,...,\widetilde{\widetilde{a}}%
_{\pi }\left( H_{m_{p}}\right) \right) \text{,}
\end{equation*}%
according respectively to (\ref{FCoef}) and (\ref{deff}).

\bigskip

\begin{theorem}
\label{T : VecCV}Suppose that, for any $j=1,...,p$, the coefficients $%
\left\{ C_{\pi }:\pi \in \hat{G}\right\} $ (as defined in (\ref{vark}))
verify either one of conditions (\ref{CLTher})-(\ref{Cond3}) (with $m_{j}$
substituting $m$). Then,%
\begin{equation}
\left\{ T\text{ \ };\text{ \ }\left( \widetilde{\widetilde{a}}_{\pi }\left(
H_{m_{1}}\right) ,...,\widetilde{\widetilde{a}}_{\pi }\left(
H_{m_{p}}\right) \right) \right\} \overset{law}{\underset{\left\{ \pi
\right\} }{\rightarrow }}\left\{ T\text{ \ };\text{ \ }\left(
N_{1}+iN_{1}^{\prime },...,N_{p}+iN_{p}^{\prime }\right) \right\}
\label{CVj}
\end{equation}%
where $\mathbf{N}_{p}=\left( N_{1},...,N_{p}\right) $ and $\mathbf{N}%
_{p}^{\prime }=\left( N_{1}^{\prime },...,N_{p}^{\prime }\right) $ are two
independent vectors of $\mathcal{N}\left( 0,1/2\right) $ i.i.d. random
variables, such that $\mathbf{N}_{p}$ and $\mathbf{N}_{p}^{\prime }$ are
independent of $T$. On the other hand, if the asymptotic relation (\ref{CVj}%
) holds, then conditions (\ref{CLTher})-(\ref{Cond3}) are necessarily
satisfied.
\end{theorem}

\bigskip

\textbf{Remark -- }The convergence relation (\ref{CVj}) is meant in the
sense of finite dimensional distributions, that is: (\ref{CVj}) is true if,
and only if, for any $k\geq 1$ and every $\left( g_{1},...,g_{k}\right) \in
G^{k}$,%
\begin{equation}
\left( T\left( g_{1}\right) ,...,T\left( g_{k}\right) ,\widetilde{\widetilde{%
a}}_{\pi }\left( H_{m_{1}}\right) ,...,\widetilde{\widetilde{a}}_{\pi
}\left( H_{m_{p}}\right) \right) \overset{law}{\underset{\left\{ \pi
\right\} }{\rightarrow }}\left( T\left( g_{1}\right) ,...,T\left(
g_{k}\right) ,N_{1}+iN_{1}^{\prime },...,N_{p}+iN_{p}^{\prime }\right)
\label{CV2}
\end{equation}

\bigskip

\begin{proof}
For some $k\geq 1$, consider vectors $\left( g_{1},...,g_{k}\right) \in G^{k}
$ and $\left( \lambda _{1},...,\lambda _{k}\right) \in \mathbb{R}^{k}$.
Then, arguments analogous to the ones adopted in the proof of Theorem \ref{T
: Main Torus} show that
\begin{equation}
\left( \sum_{i=1}^{k}\lambda _{i}T\left( g_{i}\right) ,\widetilde{\widetilde{%
a}}_{\pi }\left( H_{m_{1}}\right) ,...,\widetilde{\widetilde{a}}_{\pi
}\left( H_{m_{p}}\right) \right) \overset{law}{=}\left( I_{1}\left( \Sigma
_{i=1}^{k}\lambda _{i}h^{g_{i}}\right) ,I_{m_{1}}\left( \widetilde{%
\widetilde{h}}_{m_{1},\pi }\right) ,...,I_{m_{p}}\left( \widetilde{%
\widetilde{h}}_{m_{p},\pi }\right) \right) \text{,}  \label{vectEQ}
\end{equation}%
where the $h^{g_{i}}$'s are given by (\ref{pd0}), and the kernels $%
\widetilde{\widetilde{h}}_{m_{j},\pi }$, $j=1,...,p$, are defined in (\ref%
{deff2}). Note that the kernel $\Sigma _{i=1}^{k}\lambda _{i}h^{g_{i}}$
(which does not depend on $\pi $) is real-valued, and therefore $I_{1}\left(
\Sigma _{i=1}^{k}\lambda _{i}h^{g_{i}}\right) $ is a real-valued Gaussian
random variable. Also, by construction the following relations hold: (i) $%
\forall j=1,...,p$, $\Re \left( I_{m_{j}}\left( \widetilde{\widetilde{h}}%
_{m_{j},\pi }\right) \right) =I_{m_{j}}\left( \Re \left( \widetilde{%
\widetilde{h}}_{m_{j},\pi }\right) \right) $ and $\Im \left( I_{m_{j}}\left(
\widetilde{\widetilde{h}}_{m_{j},\pi }\right) \right) =I_{m_{i}}\left( \Im
\left( \widetilde{\widetilde{h}}_{m_{j},\pi }\right) \right) $, and
\begin{eqnarray}
\mathbb{E}\left[ I_{m_{j}}\left( \Re \left( \widetilde{\widetilde{h}}%
_{m_{j},\pi }\right) \right) I_{m_{j}}\left( \Im \left( \widetilde{%
\widetilde{h}}_{m_{j},\pi }\right) \right) \right]  &=&0  \notag \\
\mathbb{E}\left[ I_{m_{j}}\left( \Re \left( \widetilde{\widetilde{h}}%
_{m_{j},\pi }\right) \right) ^{2}\right]  &=&\mathbb{E}\left[
I_{m_{j}}\left( \Im \left( \widetilde{\widetilde{h}}_{m_{j},\pi }\right)
\right) ^{2}\right] =\frac{1}{2};  \label{aa}
\end{eqnarray}%
(ii) $\forall 1\leq k\neq j\leq p$,
\begin{eqnarray}
\mathbb{E}\left[ I_{m_{j}}\left( \Re \left( \widetilde{\widetilde{h}}%
_{m_{j},\pi }\right) \right) I_{m_{k}}\left( \Im \left( \widetilde{%
\widetilde{h}}_{m_{k},\pi }\right) \right) \right]  &=&\mathbb{E}\left[
I_{m_{j}}\left( \Re \left( \widetilde{\widetilde{h}}_{m_{j},\pi }\right)
\right) I_{m_{k}}\left( \Re \left( \widetilde{\widetilde{h}}_{m_{k},\pi
}\right) \right) \right]   \notag \\
&=&\mathbb{E}\left[ I_{m_{j}}\left( \Im \left( \widetilde{\widetilde{h}}%
_{m_{j},\pi }\right) \right) I_{m_{k}}\left( \Im \left( \widetilde{%
\widetilde{h}}_{m_{k},\pi }\right) \right) \right] =0;  \label{bb}
\end{eqnarray}%
(iii) $\forall j=1,...,p$,%
\begin{equation}
\mathbb{E}\left[ I_{m_{j}}\left( \Re \left( \widetilde{\widetilde{h}}%
_{m_{j},\pi }\right) \right) I_{1}\left( \Sigma _{i=1}^{k}\lambda
_{i}h^{g_{i}}\right) \right] =\mathbb{E}\left[ I_{m_{j}}\left( \Im \left(
\widetilde{\widetilde{h}}_{m_{j},\pi }\right) \right) I_{1}\left( \Sigma
_{i=1}^{k}\lambda _{i}h^{g_{i}}\right) \right] =0.  \label{cc}
\end{equation}%
Now suppose that either one of conditions (\ref{CLTher})-(\ref{Cond3}) hold $%
\forall m_{j}$ ($j=1,...,p$). Then, Theorem \ref{T : Main Torus} implies
that $\forall j=1,...,p$,%
\begin{equation}
\lim_{\left\{ \pi \right\} }\mathbb{E}\left[ I_{m_{j}}\left( \Re \left(
\widetilde{\widetilde{h}}_{m_{j},\pi }\right) \right) ^{4}\right]
=\lim_{\left\{ \pi \right\} }\mathbb{E}\left[ I_{m_{j}}\left( \Im \left(
\widetilde{\widetilde{h}}_{m_{j},\pi }\right) \right) ^{4}\right] =\frac{3}{4%
}\text{,}  \label{qtr}
\end{equation}%
so that Part B of Theorem \ref{T : NP2005}, together with (\ref{aa})-(\ref%
{cc}), yield that,
\begin{eqnarray}
&&\left( I_{1}\left( \Sigma _{i=1}^{k}\lambda _{i}h^{g_{i}}\right) ,\Re
\left( I_{m_{1}}\left( \widetilde{\widetilde{h}}_{m_{1},\pi }\right) \right)
,\Im \left( I_{m_{1}}\left( \widetilde{\widetilde{h}}_{m_{1},\pi }\right)
\right) ,...\right.   \notag \\
&&\text{ \ \ \ \ \ \ \ \ \ \ \ \ \ \ \ \ \ \ \ \ \ \ \ \ \ \ \ \ \ \ }\left.
...,\Re \left( I_{m_{p}}\left( \widetilde{\widetilde{h}}_{m_{p},\pi }\right)
\right) ,\Im \left( I_{m_{p}}\left( \widetilde{\widetilde{h}}_{m_{p},\pi
}\right) \right) \right)   \label{g} \\
&&\text{ \ \ \ \ \ \ \ \ \ \ \ \ \ \ \ \ \ \ \ \ \ \ \ \ \ \ \ \ \ \ \ \ \ \
\ \ \ \ \ \ \ \ }\overset{law}{\underset{\left\{ \pi \right\} }{\rightarrow }%
}\left( I_{1}\left( \Sigma _{i=1}^{k}\lambda _{i}h^{g_{i}}\right)
,N_{1},N_{1}^{\prime },...,N_{p},N_{p}^{\prime }\right) ,  \notag
\end{eqnarray}%
where the vectors $\mathbf{N}_{p}=\left( N_{1},...,N_{p}\right) $ and $%
\mathbf{N}_{p}^{\prime }=\left( N_{1}^{\prime },...,N_{p}^{\prime }\right) $
are defined in the statement of Theorem \ref{T : VecCV}. Now note that, due
to (\ref{vectEQ}), the asymptotic relation (\ref{g}) holds $\forall \left(
\lambda _{1},...,\lambda _{k}\right) $ if, and only if, (\ref{CVj}) is
verified. The proof of the first part of Theorem \ref{T : VecCV} is
therefore concluded. To prove the last part of the statement, use the
equivalence between (\ref{CVj}) and (\ref{g}) to show that (\ref{CVj})
implies that (\ref{qtr}) holds for every $j=1,...,p$. But, due to (\ref{deff}%
) and (\ref{vectEQ}), (\ref{qtr}) is equivalent to the condition: for every $%
j=1,...,p$,
\begin{equation*}
\left[ m_{j}!\widehat{C}_{\pi ,m_{j}}\right] ^{-2}\mathbb{E}\left[ \Re
\left( \widetilde{a}_{\pi }\left( H_{m_{j}}\right) \right) ^{4}\right]
\underset{\left\{ \pi \right\} }{\rightarrow }\frac{3}{4}\text{, \ \ and \ \
}\left[ m_{j}!\widehat{C}_{\pi ,m_{j}}\right] ^{-2}\mathbb{E}\left[ \Im
\left( \widetilde{a}_{\pi }\left( H_{m_{j}}\right) \right) ^{4}\right]
\underset{\left\{ \pi \right\} }{\rightarrow }\frac{3}{4}\text{,}
\end{equation*}%
so that the proof is concluded by using once again Theorem \ref{T : Main
Torus}.
\end{proof}

\bigskip

Now define $\mathfrak{A}\triangleq \left\{ a_{\pi }:\pi \in \hat{G}\right\} $%
, where the $a_{\pi }$'s are defined according to (\ref{Fcooeff}). An
immediate consequence of Theorem \ref{T : VecCV} is the following result.

\bigskip

\begin{corollary}
\label{C : vec}Fix a vector of integers $2\leq m_{1}<\cdot \cdot \cdot
<m_{p} $, and suppose that $\forall j=1,...,p$,%
\begin{equation}
\widetilde{\widetilde{a}}_{\pi }\left( H_{m_{j}}\right) \overset{law}{%
\underset{\left\{ \pi \right\} }{\rightarrow }}N+iN^{\prime }\text{,}
\label{COrCV}
\end{equation}%
where $N,N^{\prime }\sim \mathcal{N}\left( 0,1/2\right) $ are independent.
Then,
\begin{equation*}
\left( \mathfrak{A}\text{ \ ; \ }\widetilde{\widetilde{a}}_{\pi }\left(
H_{m_{1}}\right) ,...,\widetilde{\widetilde{a}}_{\pi }\left(
H_{m_{p}}\right) \right) \overset{law}{\underset{\left\{ \pi \right\} }{%
\rightarrow }}\left\{ \mathfrak{A}\text{ \ ; \ }\left( N_{1}+iN_{1}^{\prime
},...,N_{p}+iN_{p}^{\prime }\right) \right\} \text{,}
\end{equation*}%
where $\mathbf{N}_{p}=\left( N_{1},...,N_{p}\right) $ and $\mathbf{N}%
_{p}^{\prime }=\left( N_{1}^{\prime },...,N_{p}^{\prime }\right) $ are two
independent vectors of $\mathcal{N}\left( 0,1/2\right) $ i.i.d. random
variables, such that $\mathbf{N}_{p}$ and $\mathbf{N}_{p}^{\prime }$ are
independent of $\mathfrak{A}$.
\end{corollary}

\begin{proof}
Due to Theorem \ref{T : Main Torus}, (\ref{COrCV}) holds for every $%
j=1,...,p $, if, and only if, either one of conditions (\ref{M4clt})-(\ref%
{Cond3}) are verified for every $j=1,...,p$, with $m_{j}$ replacing $m$. The
conclusion is achieved by using Theorem \ref{T : VecCV}, as well as the fact
that, by (\ref{start}) and (\ref{Fcooeff}), $\sigma \left( \mathfrak{A}%
\right) =\sigma \left( T\right) $.
\end{proof}

\section{A CLT for general $F\in L_{0}^{2}\left( \mathbb{R}%
,e^{-x^{2}/2}dx\right) \label{S : genCLT}$}

We now establish a CLT such as (\ref{CLT}) for a general real-valued
function $F\in L_{0}^{2}\left( \mathbb{R},e^{-x^{2}/2}dx\right) $. Since the
sequence of normalized Hermite polynomials $\left\{ \left( m!\right)
^{-1/2}H_{m}:m\geq 0\right\} $ defined by (\ref{Her}) is an orthonormal
basis for $L_{\mathbb{R}}^{2}(\mathbb{R},$ $\left( 2\pi \right)
^{-1/2}e^{-x^{2}/2}dx)$, the function $F$ admits a unique representation of
the form
\begin{equation}
F\left( x\right) =\sum_{m=1}^{\infty }\frac{c_{m}\left( F\right) }{m!}%
H_{m}\left( x\right) \text{, \ \ }x\in \mathbb{R}\text{,}  \label{devF}
\end{equation}%
where the coefficients $c_{m}\left( F\right) $, $m=1,2...$, are such that
\begin{equation}
c_{m}\left( F\right) =\int_{\mathbb{R}}\frac{e^{-\frac{x^{2}}{2}}}{\sqrt{%
2\pi }}H_{m}\left( x\right) F\left( x\right) dx,\text{ \ \ and \ \ }\Sigma
_{m\geq 1}\frac{c_{m}\left( F\right) ^{2}}{m!}<+\infty  \label{proprCF}
\end{equation}%
(note that the sum in (\ref{devF}) starts from $m=1$ since $F$ is centered,
i.e. $F\in L_{0}^{2}\left( \mathbb{R},e^{-x^{2}/2}dx\right) $). As a
consequence, the coefficients $\widetilde{a}_{\pi }\left( F\right) $, $\pi
\in \hat{G}$, defined in (\ref{FCoef}) can be written as
\begin{equation}
\widetilde{a}_{\pi }\left( F\right) =\sum_{m=1}^{\infty }\frac{c_{m}\left(
F\right) }{m!}\int_{G}H_{m}\left( T\left( g\right) \right) \chi _{\pi
}\left( g^{-1}\right) dg=\sum_{m=0}^{\infty }\frac{c_{m}\left( F\right) }{m!}%
\widetilde{a}_{\pi }\left( H_{m}\right)  \label{dev coeff}
\end{equation}%
where the series converges in $L_{\mathbb{C}}^{2}\left( \mathbb{P}\right) $,
and the $\widetilde{a}_{\pi }\left( H_{m}\right) $'s are given by (\ref%
{FherCOEFF}). By combining Theorem \ref{T : Main Torus}\ and Theorem \ref{T
: VecCV}, from (\ref{dev coeff}) we deduce the following result.

\bigskip

\begin{theorem}
\label{T : GenCLT}For every $\pi \neq \pi _{0}$,
\begin{equation}
\mathbb{E}\left[ \left\vert \widetilde{a}_{\pi }\left( F\right) \right\vert
^{2}\right] =\sum_{m=1}^{\infty }\left( \frac{c_{m}\left( F\right) }{m!}%
\right) ^{2}\mathbb{E}\left[ \left\vert \widetilde{a}_{\pi }\left(
H_{m}\right) \right\vert ^{2}\right] =\sum_{m=1}^{\infty }\frac{c_{m}\left(
F\right) ^{2}}{m!}\widehat{C}_{\pi ,m}.  \label{VarGenF}
\end{equation}

Suppose moreover that the following relations hold

\begin{enumerate}
\item For every $m\geq 1$,
\begin{equation*}
\lim_{\left\{ \pi \right\} }\frac{m!\widehat{C}_{\pi ,m}}{\mathbb{E}\left[
\left\vert \widetilde{a}_{\pi }\left( F\right) \right\vert ^{2}\right] }%
\rightarrow \sigma _{m}^{2}\in \left( 0,+\infty \right) ;
\end{equation*}

\item $\sum_{m\geq 1}\left\{ c_{m}\left( F\right) /m!\right\} ^{2}\sigma
_{m}^{2}\triangleq \sigma ^{2}\left( F\right) <+\infty ;$

\item For every $m\geq 2$, the coefficients $\left\{ C_{\pi }:\pi \in \hat{G}%
\right\} $ given by (\ref{vark}) verify either one of conditions (\ref{Cond2}%
) and (\ref{Cond3});

\item $\lim_{p\rightarrow +\infty }\overline{\lim }_{\left\{ \pi \right\}
}\sum_{m=p+1}^{\infty }\left\{ c_{m}\left( F\right) ^{2}/m!\right\} \widehat{%
C}_{\pi ,m}=0.$
\end{enumerate}

Then,
\begin{equation*}
\widetilde{\widetilde{a}}_{\pi }\left( F\right) \triangleq \frac{\widetilde{a%
}_{\pi }\left( F\right) }{\sqrt{\mathbb{E}\left[ \left\vert \widetilde{a}%
_{\pi }\left( F\right) \right\vert ^{2}\right] }}\overset{law}{\underset{%
\left\{ \pi \right\} }{\rightarrow }}(\sigma ^{2}\left( F\right) )^{\frac{1}{%
2}}\times \left\{ N+iN^{\prime }\right\} \text{,}
\end{equation*}%
where $N,N^{\prime }\sim \mathcal{N}\left( 0,1/2\right) $ are independent
Gaussian random variables.
\end{theorem}

\begin{proof}
Fix $p\geq 1$. Assumptions 1 and 3 in the statement imply, thanks to Theorem %
\ref{T : VecCV}, that
\begin{equation*}
\frac{1}{\sqrt{\mathbb{E}\left[ \left\vert \widetilde{a}_{\pi }\left(
F\right) \right\vert ^{2}\right] }}\left( \widetilde{a}_{\pi }\left(
H_{1}\right) ,...,\widetilde{a}_{\pi }\left( H_{p}\right) \right) \overset{%
law}{\underset{\left\{ \pi \right\} }{\rightarrow }}\left( \sqrt{\sigma
_{1}^{2}}\times \left( N_{1}+iN_{1}^{\prime }\right) ,...,\sqrt{\sigma
_{p}^{2}}\times \left( N_{p}+iN_{p}^{\prime }\right) \right) \text{,}
\end{equation*}%
where $\mathbf{N}_{p}=\left( N_{1},...,N_{p}\right) $ and $\mathbf{N}%
_{p}^{\prime }=\left( N_{1}^{\prime },...,N_{p}^{\prime }\right) $ are two
independent vectors of $\mathcal{N}\left( 0,1/2\right) $ i.i.d. random
variables. In particular, it follows that%
\begin{eqnarray*}
\Phi ^{\left( p\right) }\left( \pi \right) &\triangleq &\frac{%
\sum_{m=1}^{p}\left\{ c_{m}\left( F\right) /m!\right\} \widetilde{a}_{\pi
}\left( H_{m}\right) }{\sqrt{\mathbb{E}\left[ \left\vert \widetilde{a}_{\pi
}\left( F\right) \right\vert ^{2}\right] }}\overset{law}{\underset{\left\{
\pi \right\} }{\rightarrow }}\sum_{m=1}^{p}\left\{ c_{m}\left( F\right)
/m!\right\} \times \{\sqrt{\sigma _{m}^{2}}\times \left(
N_{m}+iN_{m}^{\prime }\right) \} \\
&&\overset{law}{=}\left[ \sum_{m=1}^{p}\left\{ c_{m}\left( F\right)
/m!\right\} ^{2}\sigma _{m}^{2}\right] ^{\frac{1}{2}}\times \left\{
N_{1}+iN_{1}^{\prime }\right\} .
\end{eqnarray*}%
Now take a uniformly bounded Lipschitz function $g:\mathbb{C\mapsto R}$,
with Lipschitz coefficient equal to one. Then
\begin{eqnarray}
&&\left\vert \mathbb{E}\left[ g\left( \widetilde{\widetilde{a}}_{\pi }\left(
F\right) \right) \right] -\mathbb{E}\left[ g\left( (\sigma ^{2}\left(
F\right) )^{\frac{1}{2}}\times \left\{ N+iN^{\prime }\right\} \right) \right]
\right\vert  \notag \\
&\leq &\left\vert \mathbb{E}\left[ g\left( \widetilde{\widetilde{a}}_{\pi
}\left( F\right) \right) \right] -\mathbb{E}\left[ g\left( \Phi ^{\left(
p\right) }\left( \pi \right) \right) \right] \right\vert  \label{ddff} \\
&&+\left\vert \mathbb{E}\left[ g\left( (\sum_{m=1}^{p}\left\{ c_{m}\left(
F\right) /m!\right\} ^{2}\sigma _{m}^{2})^{\frac{1}{2}}\times \left\{
N_{1}+iN_{1}^{\prime }\right\} \right) \right] -\mathbb{E}\left[ g\left(
\Phi ^{\left( p\right) }\left( \pi \right) \right) \right] \right\vert
\notag \\
&&+\left\vert \mathbb{E}\left[ g\left( (\sum_{m=1}^{p}\left\{ c_{m}\left(
F\right) /m!\right\} ^{2}\sigma _{m}^{2})^{\frac{1}{2}}\times \left\{
N_{1}+iN_{1}^{\prime }\right\} \right) \right] -\mathbb{E}\left[ g\left(
(\sigma ^{2}\left( F\right) )^{\frac{1}{2}}\times \left\{ N+iN^{\prime
}\right\} \right) \right] \right\vert .  \notag
\end{eqnarray}%
Now recall that $\left\{ \pi \right\} $ stands for a sequence of the type $%
\left\{ \pi _{l}:l\geq 1\right\} $, and replace $\pi $ with $\pi _{l}$ in (%
\ref{ddff}). Then, by first taking the limit as $l\rightarrow +\infty $, and
then the limit as $p\rightarrow +\infty $ in the RHS of (\ref{ddff}), we
deduce from Assumptions 2 and 4 in the statement that the LHS (\ref{ddff})
converges to zero as $l\rightarrow +\infty .$ This concludes the proof.
\end{proof}

\section{The $n$-dimensional torus \label{S : Torus}}

In this last section, we focus on the case of $G$ being the $n$-dimensional
torus $\mathbb{R}^{n}/(2\pi \mathbb{Z)}^{n}$, which we parameterize as $%
(0,2\pi ]^{n}$ with addition $mod (2\pi) $ as the group operation.
In this case, the dual space $\hat{G}$ is the class of all
applications of the
type $\mathbf{\vartheta \mapsto }\exp (i\mathbf{k}^{\prime }\mathbf{%
\vartheta )}$ where $\mathbf{\vartheta }=(\vartheta _{1},...,\vartheta
_{n})\in (0,2\pi ]^{n}$ and $\mathbf{k}=\left( k_{1},...,k_{n}\right) \in
\mathbb{Z}^{n}$ (here, we identify $\hat{G}$ with the class of its
associated characters). By using the notation introduced in (\ref{Clambda1}%
)-(\ref{Clambda5}), we have also that, for every $\mathbf{k}\in \mathbb{Z}%
^{n}$,%
\begin{equation*}
\widehat{C}_{\mathbf{k},m}\triangleq \sum_{\mathbf{j}_{1}\in \mathbb{Z}%
^{n}}\cdot \cdot \cdot \sum_{\mathbf{j}_{m}\in \mathbb{Z}^{n}}\left\{ C_{%
\mathbf{j}_{1}}\cdot \cdot \cdot C_{\mathbf{j}_{m}}\right\} \mathbf{1}_{%
\mathbf{j}_{1}+\cdot \cdot \cdot +\mathbf{j}_{m}=\mathbf{k}}.
\end{equation*}%
Moreover, for any fixed $\mathbf{l}^{\ast }\in \overline{\mathbb{Z}^{n}}$,
condition (\ref{Cond3}) in the statement of Theorem (\ref{T : Main Torus})
can be rewritten as: when $\mathbf{l\rightarrow l}^{\ast }$,%
\begin{equation}
\frac{\sup_{\mathbf{j\in }\mathbb{Z}^{n}}\widehat{C}_{\mathbf{j},m-q}%
\widehat{C}_{\mathbf{j-l},q}}{\sum_{\mathbf{a\in }\mathbb{Z}^{n}}\widehat{C}%
_{\mathbf{a},m-q}\widehat{C}_{\mathbf{a-l},q}}\rightarrow 0\text{, \ \ }%
\forall q=1,...,m-1\text{. }  \label{cond4}
\end{equation}

\bigskip

\textbf{Remark }-- Condition (\ref{cond4}) bears a clear resemblance with
Lindeberg-type assumptions for the Central Limit Theorem in a martingale
difference setting, see for instance \cite{HallHey}. Indeed, in some very
simple cases (i.e. quadratic transformations of Gaussian random fields on
the 1-dimensional torus) it seems possible to derive sufficient conditions
for the CLT\ by means of martingale approximations and the extension to
complex-valued variables of convergence results for the real-valued
martingale difference sequences. However, this approach is clearly
unfeasible for general nonlinear transforms of Gaussian random fields on
higher-dimensional tori or on abstract Abelian groups.

\bigskip

As discussed in the introduction, of particular interest for physical
applications is the case where $\left\Vert \mathbf{l}\right\Vert \rightarrow
\infty $, that is, when we analyze the behavior of high-frequency
components. We discuss two examples to illustrate the application of our
results; in both cases we assume that $C_{0}=0$ to simplify the discussion.

\bigskip

\textbf{Example 1 }(\textit{Algebraic decay on the circle}) -- With this
example we show that the CLT\ fails for general Hermite transformations,
when the angular power spectrum decays algebraically. We take $n=1$ (merely
for notational simplicity) and for all $l\in \mathbb{Z}\backslash \left\{
0\right\} ,$ we assume there exist positive constants $c_{2}>c_{1}$ and $%
\alpha >1$ such that%
\begin{equation*}
c_{1}|l|^{-\alpha }\leq C_{l}\leq c_{2}|l|^{-\alpha }\text{ ;}
\end{equation*}%
of course we thus cover any model of the form $C_{l}=1/h(|l|),$ where $%
h(l)=h_{0}+h_{1}l+...+h_{p}l^{p-1}>0$ for all $l>0$ and $1/h(|l|)$ is
summable.\

We have%
\begin{eqnarray*}
\sum_{k=-\infty }^{\infty }C_{k}C_{l-k} &\leq &c_{2}\left\{
\sum_{k=1}^{\infty }\frac{1}{k^{\alpha }}\frac{1}{(k+l)^{\alpha }}%
+\sum_{k=1}^{l-1}\frac{1}{k^{\alpha }}\frac{1}{(l-k)^{\alpha }}%
+\sum_{k=l+1}^{\infty }\frac{1}{k^{\alpha }}\frac{1}{(k-l)^{\alpha }}%
\right\} \text{ .} \\
&\leq &c_{2}\left\{ \frac{2}{l^{\alpha }}\sum_{k=1}^{\infty }\frac{1}{%
k^{\alpha }}+\sum_{k=1}^{l-1}\frac{1}{k^{\alpha }}\frac{1}{(l-k)^{\alpha }}%
\right\}  \\
&\leq &c_{2}\left\{ \frac{2}{l^{\alpha }}\sum_{k=1}^{\infty }\frac{1}{%
k^{\alpha }}+2\sum_{k=1}^{[l/2]+1}\frac{1}{k^{\alpha }}\frac{1}{%
(l-k)^{\alpha }}\right\}  \\
&\leq &c_{2}\left\{ \frac{2}{l^{\alpha }}\sum_{k=1}^{\infty }\frac{1}{%
k^{\alpha }}+\frac{2}{(l/2)^{\alpha }}\sum_{k=1}^{[l/2]+1}\frac{1}{k^{\alpha
}}\right\} \leq \frac{c_{22}}{l^{\alpha }}\text{ .}
\end{eqnarray*}%
On the other hand it is immediate to see that
\begin{equation*}
\sum_{k=-\infty }^{\infty }C_{k}C_{l-k}\geq \sup_{k\in \mathbb{Z}%
}C_{k}C_{l-k}\geq C_{1}C_{l-1}\geq \frac{c_{1}^{2}}{|l-1|^{\alpha }}\geq
\frac{c_{12}}{|l|^{\alpha }}\text{ , some }c_{12}>0\text{ .}
\end{equation*}%
Arguing by induction, we have thus shown that there exist positive sequences
$c_{2q}>c_{1q},$ $q=2,3,...$ such that
\begin{equation*}
c_{1q}|l|^{-\alpha }\leq \widehat{C}_{l,q}\leq c_{2q}|l|^{-\alpha }\text{ , }
\end{equation*}%
and%
\begin{equation*}
\frac{\max_{k}\widehat{C}_{k,1}\widehat{C}_{l-k,m-1}}{\widehat{C}_{l,m}}\geq
\frac{c_{1}c_{1,m-1}}{c_{2,m}}>0\text{ for all }l\in \mathbb{Z}\backslash
\left\{ 0\right\} \text{ ,}
\end{equation*}%
whence the necessary conditions for the Central Limit Theorem (\ref{cond4})
fail for each $m\geq 2$.

\bigskip

\textbf{Remark -- }Analogous examples where the CLT fails could be easily
provided for $n>1$, by considering for instance the spectral function%
\begin{equation*}
C_{l_{1}...l_{n}}=\frac{1}{h(|l_{1}|,...,|l_{p}|)}\text{ ,}
\end{equation*}%
for $h(.,...,.)$ a multivariate polynomial which takes nonnegative values on
the positive integers. A polynomial decay of the power spectrum is common in
physical models for the large scale structure of the Universe, for instance
in the highly popular Harrison-Zeldovich model (see \cite{PeebBook}).

\bigskip

\textbf{Example 2} (\textit{Exponential decay on the circle}) -- With this
example we show that the CLT\ holds for arbitrary Hermite transformations
when the angular power spectrum decays exponentially, up to multiplicative
algebraic factors. Assume we have
\begin{equation}
c_{1}h(|l|)\exp (-\vartheta |l|)\leq C_{l}\leq c_{2}h(|l|)\exp (-\vartheta
|l|)\text{ , }l\in \mathbb{Z}\backslash \left\{ 0\right\} \text{ ,}
\label{ex11}
\end{equation}%
for strictly positive constants $\vartheta $ and $c_{2}>c_{1},$ and where $%
h(l)=h_{0}+h_{1}l+...+h_{p}l^{p-1}>0.$

Then
\begin{equation*}
\widehat{C}_{l,2}\geq \sum_{k=1}^{l-1}C_{k}C_{l-k}\geq
c_{1}\sum_{k=1}^{l-1}h(k)h(l-k)\exp (-\vartheta |l|)\geq
c_{12}|l|^{2p+1}\exp (-\vartheta |l|)\text{ ,}
\end{equation*}%
for some constant $c_{12}>0$. Iterating this argument, we obtain by induction%
\begin{equation*}
\widehat{C}_{l,q}=\sum_{k_{1}=-\infty }^{\infty }C_{k_{1}}\widehat{C}%
_{l-k_{1},q-1}\geq \sum_{k_{1}=1}^{l-1}C_{k_{1}}\widehat{C}%
_{l-k_{1},q-1}\geq c_{1q}|l|^{qp+q-1}\exp (-\vartheta |l|)\text{ , }c_{1p}>0%
\text{ .}
\end{equation*}%
On the other hand, we have also%
\begin{equation*}
\widehat{C}_{\lambda ,q}\leq c_{2q}|l|^{qp+q-1}\exp (-\vartheta |l|)\exp
(-\vartheta |l|)\text{ , some }c_{2p}>0\text{ ,}
\end{equation*}%
because%
\begin{eqnarray*}
\widehat{C}_{l,2} &=&\sum_{|k|\leq 2l}\left\{ C_{k}C_{l-k}\right\}
+\sum_{|k|>2l}\left\{ C_{k}C_{l-k}\right\} \leq c_{2}^{\prime
}|l|^{2p+1}\exp (-\vartheta |l|)+c_{2}^{\prime \prime
}C_{l}\sum_{|k|>2l}C_{k} \\
&\leq &c_{22}|l|^{2p+1}\exp (-\vartheta |l|)\text{ ,}
\end{eqnarray*}%
and then the argument is completed by induction. Hence we have%
\begin{eqnarray*}
\sup_{\lambda \in \mathbb{Z}}\widehat{C}_{\lambda ,m-q}\widehat{C}_{l\mathbf{%
-}\lambda ,q} &\leq &c_{2q}|l|^{qp+q-1}|l|^{(m-q)p+m-q-1}\exp (-\vartheta
|l|)=c_{2q}|l|^{mp+m-2}\exp (-\vartheta |l|) \\
\sum_{\mathbf{\mu \in }\mathbb{Z}}\widehat{C}_{\mu ,m-q}\widehat{C}_{l-\mu
,q} &=&\widehat{C}_{l,m}\geq c_{1q}|l|^{mp+m-1}\exp (-\vartheta |l|)\text{ ,}
\end{eqnarray*}%
whence it is immediate to see that (\ref{cond4}) follows.

\bigskip

\textbf{Remark -- }Analogous examples where a CLT of the type (\ref{CLTher})
holds for every $m\geq 2$ could be easily provided for $n>1,$ considering
for instance the spectral function%
\begin{equation*}
C_{l_{1}...l_{n}}=h(|l_{1}|,...,|l_{p}|)\exp (-\vartheta
_{1}|l_{1}|...-\vartheta _{n}|l_{n}|)\text{ ,}
\end{equation*}%
for $h(.,...,.)$ a multivariate polynomial which takes nonnegative values on
the positive integers. An exponential decay of the angular power spectrum at
very high frequencies is expected in physical models for the CMB\ random
field, due to the so-called Silk damping (or diffusion damping) effect (see
\cite{dodelson}).\newline

\end{document}